\documentclass{amsart}

\usepackage{amsmath,amssymb,amsfonts,amsthm}
\usepackage[mathscr]{eucal}
\usepackage[usenames,dvips]{color}
\usepackage[all]{xypic}

\theoremstyle{plain}
\newtheorem{theorem}{Theorem}[section]
\newtheorem*{theorem*}{Theorem}

\newtheorem{lemma}[theorem]{Lemma}
\newtheorem{corollary}[theorem]{Corollary}

\theoremstyle{definition}
\newtheorem{definition}[theorem]{Definition}

\theoremstyle{remark}
\newtheorem{remark}[theorem]{Remark}
\newtheorem{fact}[theorem]{Fact}
\newtheorem{notation}[theorem]{Notation}
\newtheorem{scenario}[theorem]{Scenario}

\newcommand{\defn}[1]{{\emph{#1}}}

\renewcommand{\a}{\alpha}
\renewcommand \b{\beta}
\renewcommand \d{\delta}
\newcommand \e{\epsilon}
\newcommand \g{\gamma}
\renewcommand \k{\kappa}
\renewcommand \l{\lambda}
\newcommand \p{\pi}
\newcommand \s{\sigma}
\newcommand\w{\omega}

\newcommand \Iff{\,\,\, \Longleftrightarrow \,\,\,}

\newcommand\force{\Vdash}
\newcommand \forces{\force}
\newcommand \B{\mathbb{B}}
\renewcommand \P{\mathbb{P}}
\newcommand{\rest}{\restriction}

\newcommand \seq[1]{{\langle{#1}\rangle}}
\newcommand \join{\vee}
\newcommand \meet{\wedge}
\newcommand \PPP{\mathcal{P}}
\newcommand \Q{{\mathbb{Q}}}
\newcommand \N{{\mathbb{N}}}

\newcommand \R{{\mathbb{R}}}

\newcommand \otp{\text{\rm otp}\,}

\newcommand \U{\mathbb U}
\renewcommand\le{\leqslant}
\renewcommand\ge{\geqslant}
\newcommand{\bool}[1]{\mathopen{{[\![}}#1\mathclose{{]\!]}}}
\newcommand \andd{\,\,\&\,\,}
\newcommand \vphi{\varphi}
\newcommand \F{\mathscr{F}}

\renewcommand \SS{\mathcal{S}}
\newcommand \RR{\mathcal{R}}

\newcommand \pure {\le_{\textup{pur}}}
\newcommand \local {\le_{\textup{loc}}}
\newcommand \NN{\mathscr{N}}

\newcommand \II{\mathcal{I}}
\newcommand \BB{\mathcal{B}}
\newcommand \onto{\twoheadrightarrow}
\newcommand \DD{\mathbb{D}}
\newcommand \compembed{\lessdot}

\renewcommand \setminus{-}

\newcommand \ups{\upsilon}
\newcommand \vas{\varsigma}
\newcommand \vpi{\varpi}
\newcommand \vro{\varrho}
\newcommand \experp{\perp_{\textup{exp}}}
\newcommand \A{\mathcal A}

\definecolor{MyDarkBlue}{rgb}{0,0,1}

\newcommand \D{\mathbb{D}}

\newcommand \On{\mathrm{On}}
\newcommand \CC{\mathscr{C}}
\newcommand \elementary{\prec}

\begin{document}

\title{Models of Real-Valued Measurability}

\author{Sakae Fuchino}
\email{fuchino@math.cs.kitami-it.ac.jp}
%\urladdr{\href{http://math.cs.kitami-it.ac.jp/~fuchino/index-e.html}}
%{math.cs.kitami-it.ac.jp/~$\sim$fuchino/index-e.html}
\address{Department of Natural Science and Mathematics \\ College of
Engineering \\  Chubu University \\ Kasugai Aichi 487-8501, Japan}

\author{Noam Greenberg}
\email{erlkoenig@nd.edu}
%\urladdr{\href{http://www.mcs.vuw.ac.nz/~greenberg}}{www.mcs.vuw.ac.nz/~$\sim$greenberg}
\address{School of Mathematics, Statistics and Computer Science\\
Victoria University of Wellington \\
Kelburn, Wellington 6005, New Zealand}

\author{Saharon Shelah}
\email{shelah@math.huji.ac.il}
%\urladdr{\href{http://shelah.logic.at}}{shelah.logic.at}
\address{The Hebrew University of Jerusalem \\ Einstein Institute of
Mathematics \\ Edmond J. Safra Campus, Givat Ram \\ Jerusalem 91904,
Israel}
\address{Department of Mathematics \\ Hill Center-Busch Campus \\ Rutgers,
The State University of New Jersey \\ 110 Frelinghuysen Road
\\ Piscataway, NJ 08854-8019 USA}

\thanks{The first author is supported by Chubu University grant S55A.
The second and third authors were supported by the United States-Israel
Binational Science Foundation (Grant no. 2002323) and NSF grant No.
NSF-DMS 0100794. Publication no. 763 in the list of Shelah's
publications.}


\maketitle

\begin{abstract}
Solovay's random-real forcing (\cite{Solovay}) is the standard way of
producing real-valued measurable cardinals. Following questions of
Fremlin, by giving a new construction, we show that there are
combinatorial, measure-theoretic properties of Solovay's model that do
not follow from the existence of real-valued measurability.
\end{abstract}

\section{Introduction}

Solovay (\cite{Solovay}) showed how to produce a real-valued measurable
cardinal by adding random reals to a ground model which contains a
measurable cardinal. (Recall that a cardinal $\k$ is \defn{real-valued
measurable} if there is an atomless, $\k$-additive measure on $\k$ that
measures all subsets of $\k$. For a survey of real-valued measurable
cardinals see Fremlin \cite{Fremlin}.)

The existence of real-valued measurable cardinals is equivalent to the
existence of a countably additive measure on the reals which measures
all sets of reals and extends Lebesgue measure (Ulam \cite{Ulam}).
However, the existence of real-valued measurable cardinals, and
particularly if the continuum is real-valued measurable, has an array of
Set Theoretic consequences reaching beyond measure theory. For example:
a real-valued measurable cardinal has the tree property (Silver
\cite{Silver}); if there is a real-valued measurable cardinal, then
there is no rapid $p$-point ultrafilter on $\N$ (Kunen); the dominating
invariant $\mathfrak d$ cannot equal a real-valued measurable cardinal
(Fremlin). And further, if the continuum is real-valued measurable then
$\Diamond_{2^{\aleph_0}}$ holds (Kunen); and for all cardinals $\l$
between $\aleph_0$ and the continuum we have $2^\l = 2^{\aleph_0}$
(Prikry \cite{Prikry}); see \cite{Fremlin}.

On the other hand, there are other properties of Solovay's model that
have not been shown to follow from the mere existence of real-valued
measurable cardinals: for example, the covering invariant for the null
ideal $\textrm{cov}(\mathcal N)$ has to equal the continuum.

Thus, Fremlin asked (\cite[P1]{Fremlin}) whether every real-valued
measurable cardinal can be obtained by Solovay's method (the precise
wording is: suppose that $\k$ is real-valued measurable; must there be
an inner model $M\subset V$ such that $\k$ is measurable in $M$ and a
random extension $M[G]\subset V$ of $M$ which contains $\PPP\k$?). The
question was answered in the negative by Gitik and Shelah
(\cite{GiSh:582}). The broader question remains: what properties of
Solovay's model follow from the particular construction, and which
properties are inherent in real-valued measurability?

In this paper we present a new construction of a real-valued measurable
cardinal and identify a combinatorial, measure-theoretic property that
differentiates between Solovay's model and the new one.

The property is the existence of what we call \emph{general sequences} -
Definition \ref{def: general sequence}. A general sequence is a sequence
which is sufficiently random as to escape all sets of measure zero.
Standard definitions of randomness are always restricted, in the sense
that the randomness has to be measured with respect to a specified
collection of null sets (from effective Martin-L\"{o}f tests to all sets
of measure zero in some ground model). Of course, we cannot simply
remove all restrictions, as no real escapes all null sets. However, we
are interested in a notion that does not restrict to a special
collection of null sets but considers them all. One way to do this is to
change the nature of the random object - here, from a real to a long
sequence of reals, and to change the nature of escaping. We remark here
that the following definition echoes (in spirit) the characterization of
(effective) Martin-L\"{o}f randomness as a string, each of whose initial
segments have high Kolmogorov complexity.

We thus introduce a notion of forcing $\Q_\k$. We show that if $\k$ is
measurable (and $2^\k = \k^+$), then in $V^{\Q_\k}$, $\k$ (which is the
continuum) is real-valued measurable (Theorem \ref{thm: new forcing
works}). We then show that in Solovay's model, the generic (random)
sequence is general (Theorem \ref{thm: the random sequence is general});
and that in the new model, no sequence is general (Theorem \ref{thm: no
general sequences}).

\subsection{Notation}

$\PPP X$ is the power set of $X$. $A\setminus B$ is set difference.
$\subset$ denotes inclusion, not necessarily proper; $\subsetneq$
denotes proper inclusion.

The reals $\R$ are identified with Cantor space $2^\w$. If $\s\in
2^{<\w}$ then $[\s] = \{ x\in \R\,:\, \s\subset x\}$ denotes the basic
open set determined by $\s$. If $\l\in \On$ then $\R^\l$ is the
$\l$-fold product of $\R$. If $\a<\l$ and $B$ is a Borel subset of $\R$
then $B^\a$ denotes $\{\bar x\in \R^\l\,:\, x_\a\in B\}$.

If $A$ is a Borel set (on some copy of Cantor space) and $W$ is an
extension of the universe $V$ then we let $A^W$ denote the
interpretation in $W$ of any code of $A$.

If $\P = (\P, \le)$ is a partial ordering then we sometimes write
$\le_\P$ for $\le$.

If $\a<\b$ are ordinals then $[\a,\b) = \{ \g\,:\, \a\le \g <\b\}$.

If $X$ and $Y$ are sets and $B\subset X\times Y$, then for $x\in X$,
$B_x = \{ y\in Y\,:\, (x,y)\in B\}$ and $B^y = \{ x\,:\, (x,y)\in B\}$
are the sections.

Suppose that $\seq{X_\a}_{\a<\d}$ is an increasing sequence of things
(ordinals, sets (under inclusion), etc.); for limit $\b\le \d$ we let
$X_{<\b}$ be the natural limit of $\seq{X_\a}_{\a<\b}$ (the supremum,
the union, etc.), and for successor $\b= \a+1$ we let $X_{<\b}= X_\a$.

\subsubsection{Forcing}

For notions of forcing, we use the notation common in the
$\textup{World} \setminus \{\textup{Jerusalem}\}$. Thus, $q\le p$ means
that $q$ extends $p$. As far as $\P$-names are concerned, we often
confuse between canonical objects and their names. Thus, $G$ is both a
generic filter but also the name of such a filter.

If $\B$ is a complete Boolean algebra and $\vphi$ is a formula in the
forcing language for $\B$, then we let $\bool{\vphi}_\B$ be the Boolean
value of $\vphi$ according to $\B$; this is the greatest element of $\B$
forcing $\vphi$. For a complete Boolean algebra the partial ordering
corresponding to $\B$ is not $\B$ itself but $\B\setminus\{0_\B\}$.
Nevertheless we often think as if the partial ordering in the forcing
were $\B$ and let $0\force_{\B}{\vphi}$ for all formulas $\vphi$ in the
forcing language.

If $\P$ is a partial ordering and $p\in \P$ then $\P(\le p)$ is the
partial ordering inherited from $\P$ on $\{q\in \P\,:\, q\le p\}$.

$\P \compembed \Q$ denotes the fact that $\P$ is a complete suborder of
$\Q$. If $\P\compembed \Q$ and $G$ is the (name for the) $\P$-generic
filter, then $\Q/G$ is the (name for the) quotient of $\Q$ by $G$: the
collection of all $q\in \Q$ which are compatible with all $p\in G$.

If $\P\subset \Q$, a strong way of getting $\P\compembed \Q$ is having a
\defn{restriction map} $q\mapsto q\rest \P$ from $\Q$ to $\P$:  a map
which is order preserving (but does not necessarily preserve $\not\le$),
and such that for all $q\in \Q$, $q\rest \P \force_{\P} q\in \Q/G$. If
$\B$ is a complete subalgebra of a complete Boolean algebra $\DD$ then
there is a restriction map from $\DD$ to $\B$; $d\rest \B = \prod^{\B}\{
b\in \B\,:\,b\ge d\}$ is in fact the largest $b\in \B$ which forces that
$d\in \DD/G$; $\DD/ G = \{ d\in \DD\,:\, d\rest \B \in G\}$.

If $\B$ is a complete subalgebra of a complete Boolean algebra $\DD$
then we let $\DD:G$ be the (name for the) quotient of $\DD$ by the
filter generated by the generic ultrafilter $G\subset \B$; $\DD:G$ is
the completion of the partial ordering $\DD/G$.

\subsection{Measure theory}

\subsubsection*{Notation; recollection of basic notions}

Recall that a \defn{measurable space} is a set $X$ together with a \defn{measure
algebra on $X$}: a countably complete Boolean subalgebra of $\PPP X$, that is some
$\SS\subset \PPP X$ containing $0$ and $X$ and closed under complementation and unions
(and intersections) of countable subsets of $\SS$. A \defn{probability measure} on a
measure space $(X,\SS)$ is a function $\mu\colon\SS\to [0,1]$ which is monotone and
countably additive: $\mu(0)=0, \mu(X)=1$ and whenever $\{B_n\,:\,n<\w\}\subset \SS$ is
a collection of pairwise disjoint sets, then $\mu(\cup B_n) = \sum \mu(B_n)$. All
measures we encounter in this work are probability measures.

Let $\mu$ be a measure on a measurable space $(X,\SS)$. Then a \defn{$\mu$-null} set
is a set $A\in \SS$ such that $\mu(A)=0$. We let $\II_\mu$ be the collection of
$\mu$-null sets; $\II_\mu$ is a countably complete ideal of the Boolean algebra $\SS$;
we can thus let $\B_\mu = \SS/ \II_\mu$; this is a complete Boolean algebra and
satisfies the countable chain condition. For $A\in \SS$, we let $[A]_\mu = A + \II_\mu
\in \B_\mu$. We often confuse $A$ and $[A]_\mu$, though. We let $\subset_\mu, =_\mu$
etc.\ be the pullback of the Boolean notions in $\B_\mu$. Namely: $A\subset_\mu B$ if
$[A]_\mu \le_{\B_\mu} [B]_\mu$ (iff $A\setminus B \in \II_\mu$), etc. We also think of
$\mu$ as measuring the algebra $\B_\mu$; we let $\mu([A]_\mu) = \mu(A)$.

\begin{definition}
Let $\SS\subset \RR$ be two measure algebras on a space $X$, and let $\mu$ be a
measure on $\SS$ and $\nu$ be a measure on $\RR$. We say that $\nu$ is
\defn{absolutely continuous with respect to $\mu$} (and write $\nu \ll \mu$) if
$\II_\mu \subset \II_\nu$; that is, if for all $A\in \SS$, if $\mu(A)=0$ then
$\nu(A)=0$. \end{definition} (Of course, if $\nu \ll \mu$, $A\in \SS$ and $\mu(A)=1$
then $\nu(A)=1$).

If $\nu \ll \mu$ then the identity $\SS\subset \RR$ induces a map $i\colon\B_\mu\to
\B_\nu$ which is a complete Boolean homomorphism. If $\II_\mu = \II_\nu \cap \SS$ then
$i$ is injective.

\begin{definition}
Let $\mu$ be a measure on $(X,\SS)$, and let $A\in \SS$ be a $\mu$-positive set. We
let $\mu \| A$, the \defn{localization} of $\mu$ to $A$, be $\mu$ restricted to $A$,
recalibrated to be a probability measure: it is the measure on $(X,\SS)$ defined by
$(\mu \| A) (B) = \mu(B \cap A)/\mu (A)$. \end{definition}

If $A=_\mu A'$ then $\mu\| A = \mu\| A'$ so we may write $\mu\| a$ for
$a\in \B_\mu$. We have $ \mu\| a \ll \mu$ and $\B_{\mu\|a} \cong
\B_\mu(\le a)$; under this identification, the natural map
$i\colon\B_\mu \to \B_{\mu\|a}$ is given by $i(b) = b\cap a$. If $a\ne
1$ (so $\mu\ne \mu\| a$) then $i$ is not injective.

\subsubsection*{Products of measures}

If for $i<2$, $\mu_i$ is a measure on a measurable space $(X_i,\SS_i)$,
then there is a unique measure $\mu_0\,\mu_1 = \mu_0\times \mu_1$
defined on the measure algebra on $X_0\times X_1$ generated by the
\defn{cylinders}, i.e.\ the sets $A_0 \times A_1$ for $A_i\in \SS_i$,
such that $(\mu_0\, \mu_1)(A_0\times A_1) = \mu(A_0)\mu_1(A_1)$ for all
cylinders $A_0\times A_1$. We recall Fubini's theorem: For any
measurable $A\subset X_0\times X_1$, we have \[ (\mu_0\,\mu_1)(A) =
\int_{X_0} \mu_1(A_x) \,d\mu_0(x),\] where for $x\in X_0$, $A_x = \{
y\in X_1 \,:\, (x,y)\in A\}$ is the $x$-section of $A$.

We note that localization commutes with finite products: \[(\mu_0 \| B_0) \, (\mu_1 \|
B_1) = (\mu_0 \, \mu_1) \| (B_0\times B_1).\]

We can generalize the notion of absolute continuity.

\begin{definition}[Generalized absolute continuity]
Suppose that $\mu$ measures $(X,\SS)$ and $\nu$ measures $(Y,\RR)$, and
further that there is a Boolean homomorphism $i\colon \SS\to \RR$. We
say that $\nu \ll \mu$ if whenever $A\in \SS$ and $\mu(A)=0$ then
$\nu(i(A)) = 0$.
\end{definition}

If $i$ is injective then we don't really get anything new (we may identify $\SS$ with
its image). In any case, the map $i$ induces a Boolean homomorphism from $\B_\mu$ to
$\B_\nu$.

The standard example is of course if $\SS= \SS_0$ and $\RR$ is the algebra generated
by $\SS_0\times \SS_1$ as above. We then let $i(A) = A\times X_1$ and get
$\mu_0\,\mu_1 \ll \mu_0$. The map  $i$ is injective and induces a complete embedding
\[i_{\mu_0}^{\mu_0\mu_1} \colon \B_{\mu_0}\to \B_{\mu_0\, \mu_1}.\]

The following is an important simplification in notation.

\begin{notation} Unless otherwise stated, we identify $\B_{\mu_0}$ with its image
under $i_{\mu_0}^{\mu_0\mu_1}$. Thus $A\in \SS_0$ is identified with $A\times X_1$.
\label{not: product identification}
\end{notation}

Thus if $A_i \in \SS_i$ then $A_0\cap A_1 = A_0\times A_1$.

\medskip

The restriction map from $\B_{\mu_0\mu_1}$ onto $\B_{\mu_0}$ is nicely
defined: for measurable $A\subset X_0\times X_1$, we let
\[ A\rest \mu_0 = \{ x\in X_0\,:\, \mu_1(A_x)>0 \}; \] this is the measure-theoretic
projection of $A$ onto $X_0$. If $A =_{\mu_0 \mu_1} A'$ then
$A\rest\mu_0 =_{\mu_0} A'\rest \mu_0$, so we indeed get a map from
$\B_{\mu_0 \mu_1}$ onto $\B_{\mu_0}$, and $[A]_{{\mu_0 \mu_1}} \rest
\B_{\mu_0} = [A\rest \mu_0]_{{\mu_0}}$.

We make use of the following.

\begin{lemma} Let $\nu$ be a measure on $X$ and for $i<2$ let $\mu_i$ be a measure
on $Y_i$. Let $B_i\subset X\times Y_i$ and let $A_i = B_i \rest \nu$.
Then $A_0\cap A_1 =_\nu 0$ iff $B_0\cap B_1 =_{\nu \mu_0 \mu_1} 0$.
\label{lem: amalgamation of measure and reductions}
\end{lemma}

\begin{proof}
To avoid confusion, in this proof we don't use the convention \ref{not:
product identification}.

Suppose that $A_0$ and $A_1$ are $\nu$-disjoint. Then $A_0\times
Y_0\times Y_1 \cap A_1\times Y_0\times Y_1 =_{\nu\mu_0\mu_1} 0$. Also,
$B_i \subset_{\nu\mu_i} A_i\times Y_i$ so $B_i\times Y_{1-i}
\subset_{\nu\mu_0\mu_1} A_i\times Y_0\times Y_1$; it follows that
$B_0\times Y_1$ and $B_1\times Y_0$ are $\nu\mu_0\mu_1$-disjoint.

Suppose that $B_0\times Y_1$ and $B_1\times Y_0$ are
$\nu\mu_0\mu_1$-disjoint. Consider $(B_0\times Y_1)\rest \nu\mu_1$; As
$B_1\times Y_0$ is a cylinder in the product $(X\times Y_1) \times Y_0$,
we have $(B_0\times Y_1)\rest \nu\mu_1 \cap B_1 =_{\nu\mu_1} 0$.
However, $(B_0\times Y_1)\rest \nu\mu_1 = A_0\times Y_1$. Now reducing
from $X\times Y_1$ to $X$ we get $A_1 = B_1\rest \nu$ is $\nu$-disjoint
from $A_0$.
\end{proof}

\subsubsection*{Infinite products}

Iterating the two-step product, we can consider products of finitely many measures.
However, we need the more intricate notion of a product of infinitely many measures.
Countable products behave much as finite products do. Let, for $n<\w$, $\mu_n$ be a
measure on a measurable space $(X_n,\SS_n)$. Again, a \defn{cylinder} is a set of the
form $\prod_{n<\w} A_n$ for $A_n\in \SS_n$. There is a unique measure $\mu_\w = \prod
\mu_i$ on the measure algebra on $\prod X_n$ generated by the cylinders such that for
a cylinder $\prod A_n$ we have $\mu_\w (\prod A_n) = \prod_{n<\w} \mu_n(A_n)$, where
the infinite product is taken as the limit of the finite products.

Localization commutes with countable products: if $A_n\in \SS_n$ is a sequence such
that $\mu_\w (\prod A_n) >0$, then $\mu_\w \| \prod A_n = \prod (\mu_n \| A_n)$. On
the other hand, note that we can have a sequence of $A_n$s such that for each $n$,
$\mu_n(A_n)>0$, but $\mu_\w (\prod A_n)=0$; in this case we can use the measure $\prod
(\mu_n \| A_n)$, but $\mu_\w \| \prod A_n$ cannot be defined.

To better understand uncountable products, we notice that a countable product can be
viewed as a direct limit of finite products. Namely, we let a \defn{finite cylinder}
be a set of the form $\prod_{n<k} A_n \times \prod_{n\ge k} X_n$ for some $k<\w$ and
$A_n\in \SS_n$ for $n\le k$. The finite cylinders are the cylinders of $\prod_{n<k}
X_n$ under the standard identification of subsets of $\prod_{n<k} X_n$ with subsets of
$\prod_{n<\w} X_n$. The measure algebra on $\prod_{n<\w} X_n$ generated by the finite
cylinders is the same as the algebra generated by the infinite cylinders. Under the
standard identifications, the finite product measures cohere and $\mu_\w$ is the
measure generated by their union.

\medskip

Let $\l>\aleph_0$ and suppose that for $\a<\l$, $\mu_{\{\a\}}$ is a measure on a
measurable space $(X_{\{\a\}},\SS_{\{\a\}})$. For $u\subset \l$ let $X_u =
\prod_{\a\in u} X_{\{\a\}}$. If $u\subset \l$ is countable, let $\SS_u$ be the measure
algebra on $X_u$ generated by the cylinders, and let $\mu_u$ be the product
$\prod_{\a\in u} \mu_{\{\a\}}$. As discussed, we can identify $\SS_u$ with an algebra
of subsets of $X_\l$ (or more generally subsets of $X_V$ for any $u\subset V\subset
\l$) by considering the measure algebra generated by cylinders \defn{with support in
$u$}: subsets of $X_V$ of the form $\prod_{\a\in u}A_\a \times X_{V\setminus u}$, for
$A_\a\in \SS_{\{\a\}}$.

For any $V\subset \l$, we let $\SS_V$ be the union of $\SS_u$ for countable $u\subset
V$ (we note that any cylinder $A\subset X_V$ has least support). The measures $\mu_u$
cohere (i.e.\ $\mu_u$ and $\mu_v$ agree on $\SS_{u\cap v}$); thus the union of the
$\mu_u$s is a measure $\mu_V$ on $\SS_V$ (this is more immediate than the countable
case because every countable subset of $\SS_V$ lies in some $\SS_u$.) As for countable
sets, we can view $\mu_V$ as measuring subsets of any $X_W$ for $V\subset W\subset
\l$. In fact, under this identification, $\SS_V$ consists of those sets of $\SS_W$
which have support in $V$, that is, sets of the form $A\times X_{W\setminus V}$ for
some $A\subset X_V$. [We note that unlike a cylinder, a set with infinite support may
not have a minimal support: consider the set of all sequences in $2^\w$ which are
eventually 0.] The measure $\mu_V$ is determined by its values on the cylinders with
finite support; for any $V\subset W\subset \l$ we have $\mu_W = \mu_V \,
\mu_{W\setminus V}$.

\subsubsection*{General framework}

For our work, we fix $\l>\aleph_0$. For all $u\subset \l$, we let $\R^u$
be the $u$-product of Cantor space. Elements of $\R^u$ are often written
as $\bar x = \seq{x_\a}_{\a\in u}$. For \emph{countable} $u\subset \l$,
we let $\SS_u$ be the collection of Borel subsets of $\R^u$, and let
$m_u$ be Lebesgue measure on $\R^u$. For countable and uncountable
$u\subset \l$, $\SS_u$ is the algebra generated by $\SS_v$ for finite
$v\subset u$ and $m_u$ is the product $\prod_{\a\in u}m_{\{\a\}}$.

The measures we shall consider will all be localizations of products of
localizations of the $m_{\{\a\}}$:
\begin{definition} Let $u\subset \l$. A \defn{pure local product measure on $u$} is a
measure on $\SS_u$ of the form $\prod_{\a\in u}(m_{\{\a\}}\| B_\a)$ for $B_\a\in
\SS_{\{\a\}}$. A \defn{local product measure on $u$} is a measure on $\SS_u$ of the
form $\nu\| B$, where $\nu$ is a pure local product measure on $u$. \label{def: local
product measure}
\end{definition}

We will mention other measures (such as a measure witnessing that a
cardinal is real-valued measurable); but when it is clear from context
that we only mention local product measures, we drop the long name and
just refer to ``measures" and ``pure measures".

If $\mu$ is a local product measure on $u$ then we let $u^\mu = u$ and call $u$ the
\defn{support} (or \defn{domain}) of $\mu$.

\subsubsection*{Topology}

We note that every $\R^u$ is also a topological space (which can be
viewed as the Tychonoff product of $\R^{\{\a\}}$ for $\a\in u$).
However, when $u$ is uncountable, then the Borel subsets of $\R^u$
properly extend $\SS_u$. This is not a concern of ours because the
completion of any local product measure measures the Borel subsets of
$\R^u$. We thus abuse terminology and when we say ``Borel" we mean a set
in $\SS_u$; so for us, every Borel set has countable support. [In some
texts, sets in $\SS_u$ are called \defn{Baire sets}. We choose not to
use this terminology to avoid confusion between measure and category.]

Recall that a measure $\mu$ which is defined on the Borel subsets of a topological
space is \defn{regular} if for all Borel $A$, $\mu(A)$ is both the infimum of $\mu(G)$
for open $G\supset A$ and the supremum of $\mu(K)$ for compact $K\subset A$. [Thus up
to $\mu$-measure 0, each Borel set is the same as a ${\mathbf\Sigma}^0_2$ (an $F_\s$)
set and as a ${\mathbf \Pi}^0_2$ (a $G_\d$) set.] Lebesgue measure is regular, and a
localization of a regular measure is also regular. Also, regularity is preserved under
products; again note that even with uncountable products, every measurable set has
countable support and so the closed sets produced by regularity have countable
support.

\begin{corollary} Every local product measure is regular.
\label{cor: measures are regular}
\end{corollary}

\subsubsection*{Random reals}

Let $\mu$ be a local product measure. Forcing with $\B_{\mu}$ is the
same as forcing with $\II_{\mu}^+  = \SS_{u^\mu} \setminus \II_{\mu}$,
ordered by inclusion. The regularity of $\mu$ shows that the closed sets
are dense in $\II_{\mu}^+$. It follows that a generic $G\subset
\II_{\mu}^+$ is determined by
\[ \{\bar r^G\} = \cap_{B\in G} B^{V[G]}. \]  We have $B\in G$ iff $\bar r^G \in
B^{V[G]}$. We have $\bar r^G \in \bigcap\{ A^{V[G]}\,:\, A\in V\text{ is
co-null}\}$; and conversely, if $W$ is an extension of $V$ and $\bar
r\in W$ lies in $\bigcap\{ A^{V[G]}\,:\, A\in V\text{ is co-null}\}$,
then $G = \{ A\in   \II^+_\mu\,:\, \bar r\in A^{W}\}$ is generic over
$V$ and $\bar r = \bar r^G$.

Suppose that $\nu,\mu$ are local product measures and that $u^\nu\cap
u^\mu=0$. Then $\nu\mu$ is a local product measure. Recall that we have
a complete embedding $i_{\nu}^{\nu\mu} \colon \B_\nu \to \B_{\nu\mu}$.
Thus if $G\subset \B_{\nu\mu}$ is generic then $G_\nu =
{(i_{\nu}^{\nu\mu})}^{-1} G$ is generic for $\B_\nu$. In fact, $\bar
r^{G_\nu} = \bar r^G \rest u^\nu$.

\subsubsection*{Quotients are measure algebras}

Let $V[G]$ be any generic extension of $V$. There is a canonical
extension of $\mu$ to a measure on $\SS_u^{V[G]}$, which we denote by
$\mu^{V[G]}$. For if $\mu = (\prod_{\a\in u}(m_{\{\a\}}\|B_\a))\|B$ then
we can let $\mu^{V[G]} = (\prod_{\a\in
u}(m_{\{\a\}}\|B_\a^{V[G]}))\|{B}^{V[G]}$. The usual absoluteness
arguments show that indeed $\mu^{V[G]}$ is an extension of $\mu$, and
does not depend on the presentation of $\mu$.

\medskip

Again let $\mu$ and $\nu$ be local product measures on disjoint
$u=u^\nu, v= u^\mu \subset \l$. We make use of the following.

\begin{fact} Let $G\subset \B_{\nu}$ be generic. Then the map $A\mapsto {A^{V[G]}}_{\bar r^G}$
induces an isomorphism from $\B_{\nu\mu}\colon G$ to $\B_{\mu^{V[G]}}$.
\end{fact}

In particular, $\force_{\B_{\nu}} ``\B_{\nu\mu}\colon G\text{ is a
measure algebra}".$

\begin{proof} Let $\p_{\nu,\mu} \colon \B_{\nu\mu} \to \B_{\nu\mu}\colon G$ be the
quotient map. We know that \[\p_{\nu,\mu}^{-1} (\B_{\nu\mu}\colon G -
\{0\}) = \B_{\nu\mu} / G \] (the partial ordering). Thus for $A \in
\SS_{u\cup v}$, we have
\begin{eqnarray*} \p_{\nu,\mu}([A]_{\nu\mu}) > 0  \Iff
    [A\rest\nu]_{\nu} \in G & \Iff & \\
        \bar r^G \in (A\rest \nu)^{V[G]} & \Iff &
    \mu^{V[G]} ({{A}^{V[G]}}_{\bar r^G}) >0,
\end{eqnarray*}
The last equivalence follows from the fact that $(A\rest \nu)^{V[G]} =
A^{V[G]} \rest \nu^{V[G]}$; again we use absoluteness. Thus we may
define an embedding $\s_{\nu,\mu}\colon \B_{\nu\mu}\colon G \to
\B_{\mu^{V[G]}}$ by letting $\s_{\nu,\mu} (\p_{\nu,\mu}([A]_{\nu\mu})) =
\left[{{A}^{V[G]}}_{\bar r^G}\right]_{\mu^{V[G]}}$. It is clear that
$\s_{\nu,\mu}$ preserves the Boolean operations.

$\s_{\nu,\mu}$ is onto: every set in the random extension is determined
by a set in the plane in the ground model (see \cite[3.1]{Judah}). For
any countable $v'\subset v$, every $\B_\nu$-name $y$ for an element of
$\R^{v'}$ corresponds to a $\mathbf{\Pi}^0_3$ function $f_y \colon
\R^u\to \R^{v'}$ defined by $f_y(x)(i)(n)=k \iff x\in
\bool{y(i)(n)=k}_{\B_\nu}$ ($f_y$ can be taken to be $\mathbf{\Pi}^0_3$
because $\bool{y(i)(n)=k}$ can be taken to be either a
$\mathbf{\Sigma}^0_2$ or a $\mathbf{\Pi}^0_2$ set.) The function $f_y$
has the property that $f_y^{V[G]}(\bar r^G) = y^G$. Let $C$ be a
$\B_\nu$-name for a Borel subset of $\R^v$. The algebra $\B_\nu$ is
c.c.c., so in $V$, there is some countable $v'\subset v$ and some
$\B_\nu$-name $C'$ for a Borel subset of $\R^{v'}$ such that
$\force_{\B_\nu} C = C' \times \R^{v\setminus v'}$. We can let
\[ A = \{ (x,f_y(x))\,:\,x\in \bool{y\in C'}_{\B_\nu}\}\,\times \,\R^{v\setminus v'} \]
where $f_y$ ranges over $\mathbf{\Pi}^0_3$ functions from $\R^u$ to
$\R^{v'}$; thus $A$ is Borel and $(A^G)_{\bar{r}^G}=C^G$. [However, in
the sequel, we do not use the fact that $\s_{\nu,\mu}$ is onto.]
\end{proof}

\subsubsection*{Commuting diagrams}

We thus have the following diagram:

\begin{equation*} \xymatrix{ \B_{\nu \mu}
\ar[rr]^{\p_{\nu,\mu}} & & \B_{\nu \mu} \colon G \ar[rr]^{\s_{\nu,\mu}} & &
\B_{\mu^{V[G]}}. & }
\end{equation*}

Suppose now that $\nu$ is a local product measure on $u$; $\mu, \vro$
are local product measures on ${v_0}, {v_1}$, and $u, v_0, v_1$ are
pairwise disjoint. Let $\ups = \mu\vro$. Let $G\subset \B_{\nu}$ be
generic. For the rest of the section, we retract our convention
\ref{not: product identification}. We thus have a complete embedding
$i_{\nu\mu}^{\nu\ups}\colon \B_{\nu\mu} \to \B_{\nu\ups}$. This
embedding induces a complete embedding $\iota_{\nu\mu}^{\nu\ups}\colon
\B_{\nu\mu}\colon G \to \B_{\nu\ups}\colon G$.

\begin{lemma} The following diagram commutes.
\[
\xymatrix{      \B_{\nu\mu}    \ar[rr]^{\p_{\nu,\mu}} \ar[dd]_{i_{\nu\mu}^{\nu\ups}} &
& \B_{\nu\mu}\colon G \ar[dd]_{\iota_{\nu\mu}^{\nu\ups}} \ar[rr]^{\s_{\nu,\mu}}
& &  \B_{\mu^{V[G]}}  \ar[dd]^{i_{\mu^{V[G]}}^{\ups^{V[G]}}}    \\
   &  &         &      \\
\B_{\nu\ups}      \ar[rr]^{\p_{\nu,\ups}}   &  & \B_{\nu\ups}\colon G
\ar[rr]^{\s_{\nu,\ups}}    &    &
\B_{\ups^{V[G]}} \\
} \]
\label{lem: diagram 1}
\end{lemma}

\begin{proof}
Let $A\in \SS_{u\cup v_0}$, and let
\begin{equation*} \begin{array}{llll}
& a & = & \p_{\nu,\mu}([A]_{\nu\mu}); \\
& a' &  = & {\iota_{\nu\mu}^{\nu\ups}} (a) = \p_{\nu,\ups}\left(\left[A\times \R^{v_1}\right]_{\nu\ups}\right); \\
& b & =  & \s_{\nu,\mu}(a) =  \left[{A^{V[G]}}_{\bar r^G}\right]_{\mu^{V[G]}};  \text{   and} \\
& b' & =  & \s_{\nu,\ups}(a') = \left[{{(A\times
\R^{v_1})}^{V[G]}}_{\bar r^{G}}\right]_{\ups^{V[G]}}.
\end{array}
\end{equation*}
The desired equation ${i_{\mu^{V[G]}}^{\ups^{V[G]}}}(b) = b'$ follows
from the fact that \[{{A}^{V[G]}}_{\bar r^G} \times {\R^{v_1}}^{V[G]} =
{{\left(A\times \R^{v_1}\right)}^{V[G]}}_{\bar r^{G}}. \qedhere\]
\end{proof}

Note that ${i_{\mu^{V[G]}}^{\ups^{V[G]}}}$ is measure-preserving.

\

Next, suppose that $\vas, \vro$ are local product measures on $u_0, u_1$ and that
$\mu$ is a local product measure on $v$; and that $u_0,u_1,v$ are pairwise disjoint.
We let $\nu = \vas\vro$.

As $i_{\vas}^{\nu}$ is a complete embedding, we know that if $G_\nu \subset \B_\nu$ is
generic, then $G_\vas = {(i_{\vas}^{\nu})}^{-1} G_\nu$ is also generic. The map
$i_{\vas\mu}^{\nu\mu}$ induces a complete embedding $\iota_{\vas\mu}^{\nu\mu}$ from
$\B_{\vas\mu}\colon G_\vas$ to $\B_{\nu\mu}\colon G_\nu$.

Also, as $V[G_\vas] \subset V[G_\nu]$ we have (relying on absoluteness) a
\emph{measure-preserving} embedding $\i_{\vas,\mu}^{\nu,\mu}\colon
\B_{\mu^{V[G_\vas]}} \to \B_{\mu^{V[G_\nu]}}$, given by $[B]_{\mu^{V[G_\vas]}} \to
\left[B^{V[G_\nu]}\right]_{\mu^{V[G_\nu]}} $.

\begin{lemma}
The following diagram commutes:
\begin{equation*}
\xymatrix{      \B_{\vas\mu}    \ar[rr]^{\p_{\vas,\mu}} \ar[dd]_{i_{\vas\mu}^{\nu\mu}}
& & \B_{\vas\mu}\colon G_\vas \ar[dd]_{\iota_{\vas\mu}^{\nu\mu}}
\ar[rr]^{\s_{\vas,\mu}} & & \B_{\mu^{V[G_\vas]}}
\ar[dd]_{\i_{\vas,\mu}^{\nu,\mu}}  \\
   &  &         &      \\
\B_{\nu\mu}      \ar[rr]^{\p_{\nu,\mu}}   &  &
\B_{\nu\mu}\colon G_\nu \ar[rr]_{\s_{\nu,\mu}}    &  & \B_{\mu^{V[G_\nu]}}  \\
}
\end{equation*}
\label{lem: diagram 2}
\end{lemma}

\begin{proof}
Let $A\in \SS_{u_0\cup v}$. We let:
\begin{equation*} \begin{array}{llll}
& a & = & \p_{\vas,\mu}\left([A]_{\vas\mu}\right); \\
& a' &  = & {\iota_{\vas\mu}^{\nu\mu}} (a) = \p_{\nu,\mu}\left(\left[A\times \R^{u_1}\right]_{\nu\mu}\right); \\
& b & =  & \s_{\vas,\mu}(a) =  \left[{A^{V[G_\vas]}}_{\bar r^{G_\vas}}\right]_{\mu^{V[G_\vas]}};  \text{   and} \\
& b' & =  & \s_{\nu,\mu}(a') = \left[{{\left(A\times
\R^{u_1}\right)}^{V[G_\nu]}}_{\bar r^{G_\nu}}\right]_{\mu^{V[G_\nu]}}.
\end{array}
\end{equation*}
We want to show that ${\i_{\vas,\mu}^{\nu,\mu}}(b) = b'$. Letting $B=
{A^{V[G_\vas]}}_{\bar r^{G_\vas}}$ and $B' = {{\left(A\times
\R^{u_1}\right)}^{V[G_\nu]}}_{\bar r^{G_\nu}}$, we show that $B' =
B^{V[G_\nu]}$. We know, though, that $\bar r^{G_\nu} = {{\bar
r}^{G_\vas}} {}^{\smallfrown} {{\bar r}^{G_\vro}}$, from which we deduce
that $B' = {A^{V[G_\nu]}}_{\bar r^{G_\vas}}$. The conclusion follows
from absoluteness.
\end{proof}

In our third scenario, we have $\vro, \mu$ which are local product measures on
disjoint $v,u$; and we let $\nu = \vro\| B$ be some localization of $\vro$. In this
case we have a projection $i_{\vro}^{\nu}\colon \B_\vro \onto \B_\nu$. Let
$G_\nu\subset \B_\nu$ be generic; then $G_\vro = {(i_\vro^\nu)}^{-1} G_\nu$ is
generic, but in fact contains no less information; so we denote the extension by
$V[G]$.

\begin{lemma}
The following diagram commutes:
\begin{equation*}
\xymatrix{ \B_{\vro\mu}    \ar[rr]^{\p_{\vro,\mu}} \ar[dd]_{i_{\vro\mu}^{\nu\mu}} & &
\B_{\vro\mu}\colon G_\vro \ar[dd]_{\iota_{\vro\mu}^{\nu\mu}} \ar[drr]^{\s_{\vro,\mu}}
& & \\
&  & &  &
\B_{\mu^{V[G]}}   \\
\B_{\nu\mu}      \ar[rr]^{\p_{\nu,\mu}} &  & \B_{\nu\mu}\colon G_\nu
\ar[urr]_{\s_{\nu,\mu}} &  & }
\end{equation*}
\label{lem: diagram 3}
\end{lemma}

\begin{proof}
As usual we take $A\in \SS_{v \cup u}$ and follow $[A]_{\vro\mu}$ along
the diagram. We have
$\iota_{\vro,\mu}^{\nu,\mu}(\p_{\vro,\mu}([A]_{\vro\mu})) =
\p_{\nu,\mu}([A]_{\nu\mu})$; and
$\s_{\vro,\mu}(\p_{\vro,\mu}([A]_{\vro\mu})) = \left[{A^{V[G]}}_{\bar
r^{G_{\vro}}}\right]_{\mu^{V[G]}}$ and
$\s_{\nu,\mu}(\p_{\nu,\mu}([A]_{\nu\mu})) = \left[{A^{V[G]}}_{\bar
r^{G_{\nu}}}\right]_{\mu^{V[G]}}$; the latter two are equal because
$\bar r^{G_\vro} = \bar r^{G_\nu}$.
\end{proof}

Our last case is perhaps the easiest (in fact we do not use it later but
we include it for completeness.) Suppose that $u$ and $v$ are disjoint
and that $\nu,\mu$ are local product measures on $u,v$ respectively.
Suppose that $C\in \BB_{\mu}$ is positive; let $\ups = \mu\| C$. Let
$G\subset \B_\nu$ be generic. Then by absoluteness $\ups^{V[G]} =
\mu^{V[G]} \| C^{V[G]}$. Note that unlike the previous cases, the
Boolean homomorphism $i_{\mu^{V[G]}}^{\ups^{V[G]}}$ is not
measure-preserving.

\begin{lemma} The following diagram commutes.
\[
\xymatrix{      \B_{\nu\mu}    \ar[rr]^{\p_{\nu,\mu}} \ar[dd]_{i_{\nu\mu}^{\nu\ups}} &
& \B_{\nu\mu}\colon G \ar[dd]_{\iota_{\nu\mu}^{\nu\ups}} \ar[rr]^{\s_{\nu,\mu}}
& &  \B_{\mu^{V[G]}}  \ar[dd]^{i_{\mu^{V[G]}}^{\ups^{V[G]}}}    \\
   &  &         &      \\
\B_{\nu\ups}      \ar[rr]^{\p_{\nu,\ups}}   &  & \B_{\nu\ups}\colon G
\ar[rr]^{\s_{\nu,\ups}}    &    &
\B_{\ups^{V[G]}} \\
} \]
\label{lem: diagram 4}
\end{lemma}

\begin{proof}
Immediate, because for $A\in \SS_{u\cup v}$,
$i_{\nu\mu}^{\nu\vas}([A]_{\nu\mu}) = [A]_{\nu\ups}$ (which is the same
as $[A\cap (\R^u\times C)]_{\nu\ups}$).
\end{proof}

%%%%%%%%%%%%%%%%%%%%%%%%%%%%%%%%%%%%%%%%%%%%%%%%%%%%%%%%%%%%%%%%%%%%%%%%%%%%%%%%%%%%

\section{Solovay's construction}

We hope that the gentle reader will not be offended if we repeat a proof
of Solovay's original construction of a real-valued measurable cardinal,
starting from a measurable cardinal. The exposition which we give is
different from the one found in most textbooks, indeed from the one
given by Solovay in his paper; since in the rest of this paper we shall
elaborate on this proof, we thought such an exposition may be useful.

Let $\k$ be a measurable cardinal; let $j\colon V\to M$ be an elementary
embedding of $V$ into a transitive class model $M$ with critical point
$\k$, such that $M^\k \subset M$.

We move swiftly between $M$, $V$, $M[G]$ and $V[G]$. Whenever necessary
we indicate where we work, but many notions are absolute and there is
not much danger of confusion.

The forcing Solovay uses is $\P=\B_{m_\k}$, i.e.\ forcing with Borel
subsets of $\R^\k$ of positive Lebesgue measure. We show that after
forcing with $\P$, $\k$ is real-valued measurable.

We have $j(\P) = \left(\B_{m_{j(\k)}}\right)^M = \B_{m_{j(\k)}}$. Also,
$\P\in M$ and $ \P = \left(\B_{m_\k}\right)^M$.

Let $G\subset \P$ be generic over $V$. Then $G$ is generic over $M$. We
have the following diagram:
\begin{equation*} \xymatrix{ j(\P) \ar[rrr]^{\p_{m_\k,
m_{[\k,j(\k))}}} & & & j(\P) \colon G \ar[rrr]^{\s_{m_\k,
m_{[\k,j(\k))}}} & & & \left(\B_{{m_{[\k,j(\k))}}}\right)^{M[G]}.
& }
\end{equation*}
For shorthand, we let $\p = \p_{m_\k, m_{[\k,j(\k))}}$ and we let $\nu$
be the pullback to $j(\P)\colon G$ of ${m_{[\k,j(\k))}}^{M[G]}$ by
${\s_{m_\k, m_{[\k,j(\k))}}}$.

\medskip

Let $A$ be a $\P$-name for a subset of $\k$. In $M$, $j(A)$ is a
$j(\P)$-name for a subset of $j(\k)$. Let $b_A = \left(\bool{\k\in
j(A)}_{j(\P)}\right)^{M}$ (note $A\mapsto b_A$ is in $V$) and in $V[G]$
let $\mu(A) = \nu(\p(b_A))$. We now work in $V$ so we refer to the
objects defined as names.

\begin{lemma} Suppose that $a\in \P$, that $A,B$ are $\P$-names for subsets of $\k$, and
that $a\forces_\P A\subset B$. Then $a\forces_\P \mu(A) \le \mu(B)$.
\end{lemma}
\begin{proof}
The point is that $j(a)=a$. Let $G\subset \P$ be generic over $V$ such
that $a\in G$. In $M$, $a\force_{j(\P)} j(A)\subset j(B)$. Let $b = b_A
\cap a$. As $a\in G$ we have $\p(a) = 1_{j(\P)\colon G}$ so $\p(b) =
\p(b_A)$. However, in $M$, $b \force_{j(\P)} \k\in j(B)$ (as it forces
that $j(A)\subset j(B)$ and that $\k\in j(A)$) and so $b\le b_B$. It
follows that $\p(b_A) \le \p(b_B)$ so $\mu(A)\le \mu(B)$. As $G$ was
arbitrary, $a\forces_\P \mu(A)\le \mu(B)$.
\end{proof}

It follows that $\mu$, rather than being defined on names for subsets of
$\k$, can be well-defined on subsets of $\k$ in $V[G]$. The following
lemmas ensure that $\mu$ is indeed a (non-trivial) $\k$-complete
measure.

\begin{lemma} Let $A\subset \k$ be in $V$ and let $G\subset \P$ be generic over $V$.
If $\k\in j(A)$ then $\mu(A)=1$ and if $\k\notin j(A)$ then $\mu(A)=0$.
\end{lemma}

\begin{proof}
Suppose that $\k\in j(A)$. Then in $M$, $1_{j(\P)}\forces \k\in j(A)$.
Thus $b_A = 1_{j(\P)}$ so $\p(b_A) = 1_{j(\P)\colon G}$. Thus
$\mu(A)=1$. On the other hand, if $\k\notin j(A)$ then in $M$, no $b\in
j(\P)$ forces that $\k\in j(A)$, so $b_A=0$; it follows that $\mu(A)=0$.
\end{proof}

\begin{lemma} Suppose that $\seq{A_n}_{n<\w}$ is a sequence of
$\P$-names for subsets of $\k$. Suppose that $a\in \P$ forces that $A=
\bigcup_{n<\w} A_n$ is a disjoint union. Then $a\forces_\P \mu(A) =
\sum_{n<\w} \mu(A_n)$.
\end{lemma}

\begin{proof}
We have $j(a)=a$ and $j(\seq{A_n}_{n<\w}) = \seq{j(A_n)}_{n<\w}$; so in
$M$, $a$ forces (in $j(\P)$) that $j(A) = \cup_{n<\w} j(A_n)$ is a
disjoint union. Again let $G$ be generic such that $a\in G$.

Let $l,k<\w$ and $l\ne k$. In $M$, $a\forces_{j(\P)} j(A_k) \cap j(A_l)
=0$ so $b_{A_k} \cap a$ and $b_{A_l}\cap a$ are disjoint in $j(\P)$. As
$a\in G$ it follows that in $j(\P)\colon G$, $\p(b_{A_k}) \meet
\p(b_{A_l}) = 0$. We thus have $\nu (\sum^{j(\P)\colon G}_n \p(b_{A_n}))
= \sum_n \mu (A_n)$. It thus suffices to show that $\sum^{j(\P)\colon
G}_n \p(b_{A_n}) = \p(b_A)$.

For any $n<\w$, $a\forces A_n\subset A$ so as we saw before,
$\p(b_{A_n}) \le \p(b_A)$. To show the other inclusion, let $b = b_A
\setminus \sum^{j(\P)}_n b_{A_n}$. Then in $M$, $b\force_{j(\P)} \k\in
j(A)\setminus \cup_n j(A_n)$. Since $a\force_{j(\P)} j(A) \subset \cup_n
j(A_n)$ we must have $a\cap b =0$, which implies that $\p(b)=0$. The
equality follows.
\end{proof}

\begin{lemma} Suppose that $\g<\k$ and that $\seq{A_\a}_{\a<\g}$ is a sequence of
$\P$-names for subsets of $\k$. Suppose that $a\in \P$ and $a\forces_\P
\forall \a<\g\,(\mu(A_\a)=0)$. Then $a\forces_\P \mu(\cup_\a A_\a)=0$.
\end{lemma}
\begin{proof}
Let $A$ be a $\P$-name for a subset of $\k$ such that $a\forces_\P A =
\cup_{\a<\g} A_\a$. Then in $M$, $a\forces_{j(\P)} j(A) = \cup_{\a<\g}
j(A_\a)$ and also $a\forces_{j(\P)} \forall \a<\g\,\,
\left[\p(b_{A_\a})=0\right]$, that is, $a\forces \forall \a<\g\,\,
[\k\notin j(A_\a)]$. Thus $a\forces_{j(\P)} \k\notin j(A)$ so $a\cap b_A
=_{j(\P)} 0$ so $a\force_{j(\P)} \p(b_A)=0$ so $a\force_\P \mu(A)=0$.
\end{proof}

%%%%%%%%%%%%%%%%%%%%%%%%%%%%%%%%%%%%%%%%%%%%%%%%%%%%%%%%%%%%%%%%%%%%%%%%%%%

\section{A new construction of a real-valued measurable cardinal}

\begin{definition} A set of ordinals $u$ is of \defn{Easton type} if whenever $\theta$
is an inaccessible cardinal, $u\cap \theta$ is bounded below $\theta$.
\end{definition}

Let $\Q$ consist of the collection of all local product measures
(definition \ref{def: local product measure}) whose support is of Easton
type. If $u$ is a set of ordinals, then we let $\Q_u$ be the collection
of measures in $\Q$ whose support is contained in $u$.

Let $\mu,\nu\in \Q$. We say that $\nu$ is a \defn{pure} extension of
$\mu$ (and write $\nu \pure \mu$) if $\nu = \mu\, \vas$ where $\vas$ is
a \emph{pure} local product measure.

We say that $\nu$ is a \defn{local} extension of $\mu$ (and write
$\nu\local \mu$) if $\nu$ is a localization of $\mu$ (in particular
$\mu$ and $\nu$ have same support $u$).

We let $\nu$ extend $\mu$ ($\nu \le \mu$) if there is some $\vas$ such
that $\nu\local \vas\pure \mu$. It is not hard to verify that $\le$ is
indeed a partial ordering on $\Q$, and in fact on every $\Q_u$.

\begin{lemma} Suppose that $\nu \pure \vas$ and $\vas \local \mu$. Then $\nu\le \mu$.
\end{lemma}

\begin{proof}
Let $\ups$ be a pure  measure such that $\nu = \vas \,\ups$. Let $B\in
\SS_{u^{\mu}}$ such that $\vas = \mu\|B$.  Then $\nu = (\mu\, \ups)\|
B$, so $\mu\, \ups$ witnesses $\nu\le \mu$.
\end{proof}

Note that if $\nu \le \mu$ then $\nu \ll \mu$.

\subsection{Characterization of a generic}

We wish to find some characterization of a generic filter of $\Q_u$,
analogous to the description of a generic for random forcing in terms of
a random real. We need to discuss compatibility in $\Q$.

\subsubsection{Compatibility in $\Q$}

\begin{definition} Let $\mu, \nu\in \Q$. We say that $\mu$ and
$\nu$ are \defn{explicitly incompatible} (and write $\mu\experp \nu$) if
there is some $B\in \SS_{u^\mu\cap u^\nu}$ such that $\mu(B)=0$ but
$\nu(B)=1$.
\end{definition}

It is clear that if $\mu \experp \nu$ then $\mu\perp \nu$ (in $\Q$ and
in every $\Q_u$); because we cannot have some $\vas \ll \mu, \nu$.

\begin{lemma} Suppose that $u^\mu\cap u^\nu\ne 0$ and $\mu \not\experp \nu$. Then
there is some pure  measure $\vas$ on $u^\mu\cap u^\nu$ such that
$\mu,\nu\le \vas$.
\label{lem: the explicit perp}
\end{lemma}

\begin{proof}
Let $\mu_0,\nu_0$ be pure local product measures such that $\mu = \mu_0
\| C^\mu$, $\nu = \nu_0 \| C^\nu$ for some positive sets $C^\mu,C^\nu$.
Pick sequences $\seq{B_\a^\mu}_{\a\in u^\mu}, \seq{B_\a^\nu}_{\a\in
u^\nu}$ which define $\mu_0,\nu_0$ (i.e.\ $\mu_0 = \prod_{\a\in u^\mu}
(m_{\{\a\}}\|B_\a^\mu)$ and similarly for $\nu_0$). Note that $\mu_0$
and $\nu_0$ are not by any means unique, but that the $B_\a$s are
determined (up to Lebesgue measure) by $\mu_0,\nu_0$.

Let $v= u^\mu\cap u^\nu$. Let $\vas = \prod_{\a\in v} (m_{\{\a\}} \|
(B_\a^\mu \cup B_\a^\nu))$.

First we show that for all but countably many $\a\in v$ we have
$B_\a^\mu =_{m_{\{\a\}}} B_\a^\nu$. Suppose not; then for some $\e<1$ we
have some countable, infinite $w\subset v$ such that for all $\a\in w$,
$(m_{\{\a\}}\| B_\a^\mu) (B_\a^\nu) <\e$ (or the other way round). Let $
A = \prod_{\a\in w} B_\a^\nu$. Then $\nu(A)=1$ but $\mu(A)=0$.

Let $w = \{\a\in v\,:\, B_\a^\mu \ne_{m_{\{\a\}}} B_\a^\nu\}$. We assume
that $w\ne 0$ for otherwise we're done. Let $A^\mu = \prod_{\a\in w}
B_\a^\mu$ (and similarly define $A^\nu$). It is sufficient to show that
$\vas(A^\mu), \vas(A^\nu)>0$; it will then follow that $\mu_0$ is a pure
extension of $\vas\| A^\mu$, and similarly for $\nu_0$. Suppose that
$\vas( A^\mu)=0$. Let $\a\in w$; let $a_\a = m_{\{\a\}}(B_\a^\nu
\setminus B_\a^\mu)$, $c_\a = m_{\{\a\}}(B_\a^\mu \setminus B_\a^\nu)$
and $b_\a = m_{\{\a\}}(B_\a^\mu\cap B_\a^\nu)$. The assumption is that
$\prod_{\a\in w} \frac{b_\a+c_\a}{a_\a + b_\a + c_\a} =0$. However, for
each $\a\in w$, $\frac{b_\a}{a_\a+ b_\a} \le \frac{b_\a+ c_\a}{a_\a+
b_\a+ c_\a}$, which means that $\prod_{\a\in w} (m_{\{\a\}}\|
B_\a^\nu)(B_\a^\mu) =0$, so $\nu(A^\mu)=0$ (and of course
$\mu(A^\mu)=1$).
\end{proof}

For $\mu,\nu\in \Q$, if $u^\mu \cap u^\nu = 0$ then $\mu\nu \le \mu,\nu$
and so $\mu$ and $\nu$ are compatible. The following is the
generalization we need:

\begin{lemma} Let $u$ be a set of ordinals and let
$\mu,\nu\in \Q_u$. Then $\mu\perp_{\Q_u} \nu$ iff $\mu \experp \nu$.
\label{lem: perpendicularity}
\end{lemma}
\begin{proof}
Suppose that $\mu\not\experp \nu$. We may assume that $v= u^\mu \cap
u^\nu \ne 0$. By lemma \ref{lem: the explicit perp}, find some pure
$\vas$ on $v$, some pure $\mu_1,\nu_1$ and some $C^\mu, C^\nu$ such that
$\mu = (\vas \, \mu_1)\| C^\mu$, $\nu = (\vas \, \nu_1) \| C^\nu$. Let
$\ups = \vas\,\mu_1\, \nu_1$. We have $\ups(C^\mu \cap C^\nu)>0$ for
otherwise, by lemma \ref{lem: amalgamation of measure and reductions},
$C^\mu\rest \vas, C^\nu\rest \vas$ are $\vas$-disjoint and would witness
that $\mu\experp \nu$. Then $\ups \| (C^\mu \cap C^\nu)$ is a common
extension of $\mu$ and $\nu$; for example, $\ups \| (C^\mu \cap C^\nu) =
(\mu \, \nu_1)\| C^\nu$.
\end{proof}

\begin{remark} If $\mu\not\perp \nu$ then there is some $\ups$ on $u^\mu\cup u^\nu$
which is a common extension of $\mu$ and $\nu$. In fact, the common
extension constructed in the proof of lemma \ref{lem: perpendicularity}
is the greatest common extension of $\mu$ and $\nu$ in $\Q$ (thus this
extension does not depend on the choice of $\vas$).
\end{remark}

\subsubsection{Characterization of the generic}

Let $u$ be a set of ordinals, and let $G\subset \Q_u$ be generic over
$V$. Let
\[ A_G = \cap\, \{ B^{V[G]}\,:\, \textup{for some }\mu\in G, \,\,\mu(B)=1\}. \]

\begin{lemma} $A_G$ is not empty. \end{lemma}
\begin{proof}
Let $\F_G = \{ B^{V[G]}\,:\, B \textup{ is closed and for some }\mu\in
G, \,\,\mu(B)=1\}$, and let $B_G = \cap \F_G$. We show that $B_G = A_G$
and that $B_G$ is not empty.

For the first assertion, recall (corollary \ref{cor: measures are
regular}) that every $\mu\in \Q$ is a regular measure. Let $\mu \in
\Q_u$ and let $B$ be of $\mu$-measure 1. There is some closed $A\subset
B$ of positive measure, so $\mu\| A \in \Q_u$. Thus by genericity, for
every $B$ such that $\mu(B)=1$ for some $\mu\in G$, there is some closed
$A\subset B$ and some $\nu\in G$ such that $\nu(A)=1$. This shows that
$B_G = A_G$.

Next, we note that $\F_G$ has the finite intersection property. Let
$F\subset \F_G$ be finite. For $B\in F$ let $\nu_B\in G$ witness $B\in
\F_G$. There is some $\mu\in G$ which extends all $\nu_B$ for $B\in F$.
Then $\mu (\cap F) =1$ which implies that $\cap F \ne 0$. As
${\R^u}^{V[G]}$ is compact, $B_G\ne 0$.
\end{proof}

In fact,
\begin{lemma} $A_G$ is a singleton $\{\bar s^G\}$. \end{lemma}
\begin{proof}
Let $\a< \l$ and let $n<\w$. There is some $\mu\in G$ and some $\s\in
2^n$ such that $\mu ([\s]^\a) =1$. For given any $\mu$ we can extend it
to some $\nu$ such that $\a \in u^{\nu}$ and then extend $\nu$ locally
to some $\vas$ such that $\vas ([\s]^\a)=1$ for some $\s\in 2^n$.
\end{proof}

As usual,
\begin{lemma} $V[G] = V[\bar s^G]$. \end{lemma}
\begin{proof}
In fact, $G$ can be recovered from $\bar s^G$ because for all $\mu\in
\Q_u$, $\mu\in G$ iff for all $B$ such that $\mu(B)=1$ we have $\bar s^G
\in B^{V[G]}$. For if $\mu\notin G$ then there is some $\nu\in G$ such
that $\nu\perp \mu$. By lemma \ref{lem: perpendicularity}, there is some
$B$ such that $\nu(B)=0$ and $\mu(B)=1$. Then $\bar s^G\notin B^{V[G]}$.
\end{proof}

\subsubsection{The size of the continuum}

Here is an immediate application:

\begin{lemma} $\Q_u$ adds at least $|u|$ reals. \end{lemma}
\begin{proof}
Let $G$ be generic and let $\bar s^G$ be the generic sequence. We want
to show that for distinct $\a,\b\in u$ we have $\bar s^G_\a\ne \bar
s^G_\b$. Let $\mu\in G$ be such that $\a,\b\in u^\mu$. $\mu \le
m_{\{\a,\b\}}$ so $\mu \ll m_{\{\a,\b\}}$. Let $ A = \{(x,y)\in
\R^{\{\a\}}\times \R^{\{\b\}}\,:\, x\ne y\}$ be the complement of the
diagonal. Then $m_{\{\a,\b\}} (A)=1$ so $\mu(A)=1$. Thus $\mu \forces
\bar s^G\in A^{V[G]}$. But $A^{V[G]} = \{(x,y)\in
{\R^{\{\a\}}}^{V[G]}\times {\R^{\{\b\}}}^{V[G]}\,:\, x\ne y\}$. Thus
$\bar s^G_\a\ne \bar s^G_\b$.
\end{proof}

\subsection{More on local and pure extensions}

Let $\mu\in \Q$. The collection of local extensions of $\mu$ (ordered by
$\le$) is isomorphic to $\B_\mu$, so we identify the two.

\begin{lemma} Let $\mu\in \Q_u$. Then $\B_\mu \compembed \Q_u(\le \mu)$.
\label{lem: random embeds}
\end{lemma}
\begin{proof}
Let $A,B\in \B_\mu$. Then $A$ and $B$ are compatible in $\B_\mu$ iff
$\mu(A\cap B)>0$ iff $\mu\|A, \mu\|B$ are compatible in $\Q$.

Let $\seq{A_n}_{n<\w}$ be a maximal antichain of $\B_\mu$. Let $\nu \in
\Q_u$, $\nu\le \mu$. Since $\mu(\cup A_n)=1$ we have $\nu(\cup A_n)=1$
and so for some $n<\w$ we have $\nu(A_n)>0$. Then $\nu\| A_n$ is a
common extension of $\nu$ and $\mu\| A_n$.
\end{proof}

It follows that $\bar s^G$ is a string of random reals.

\begin{remark} For all $u\subset v$ we have $\Q_u \compembed \Q_v$; we do not need
this fact. \end{remark}

\begin{definition} Let $\mu\in \Q_u$ and let $\U\subset \Q_u$. We say that $\mu$
\defn{determines} $\U$ if $\U\cap \B_\mu$ is dense in $\B_\mu$.
\end{definition}

We say that $\mu \in \Q_u$ \defn{determines} a formula $\vphi$ of the
forcing language for $\Q_u$ if $\mu$ determines $\{ \nu\in \Q_u\,:\,
\nu\text{ decides }\vphi\}$. Of course, this depends on $u$, so if not
clear from context we will say ``$u$-determines". Informally, $\mu$
determining $\vphi$ means that $\vphi$ is transformed to be a statement
in the random forcing $\B_\mu$, which is a simple notion, compared to
formulas of $\Q_u$. If $\mu$ determines pertinent facts about a
$\Q_u$-name then that name essentially becomes a $\B_\mu$-name.

\medskip

For a formula $\vphi$ of the forcing language for $\Q_u$ and $\mu\in
\Q_u$ we let
\[ \bool{\vphi}^{u}_\mu  = {\textstyle \sum^{\B_\mu}} \{ b \in \B_\mu\,:\, \mu\|b
\force_{\Q_u} \vphi\}.\] Then $\mu$ $u$-determines $\vphi$ iff
$\bool{\vphi}^{u}_\mu \join \bool{\lnot\vphi}^{u}_\mu = 1_{\B_\mu}$.
Recall that if $\nu \le \mu$ then $\nu \ll \mu$ so there is a natural
map $i_{\mu}^{\nu}\colon \B_{\mu}\to \B_\nu$ (which is a
measure-preserving embedding if $\nu$ is a pure extension of $\mu$). For
all $a\in \B_\mu$, if $i_\mu^\nu(a)\ne 0$ then $\nu \| i_\mu^\nu(a)
\le_{\Q} \mu \| a$, so for all $\vphi$, $\bool{\vphi}^{u}_\nu
\ge_{\B_\nu} i_\mu^\nu(\bool{\vphi}^{u}_\mu)$. Thus if $\mu$ determines
$\vphi$ then so does $\nu$ and in this case $\bool{\vphi}^{u}_\nu =
i_\mu^\nu (\bool{\vphi}^{u}_\mu)$. If also $\nu\pure \mu$ then these
Boolean values have the same measure: $\mu(\bool{\vphi}^{u}_\mu) = \nu
(\bool{\vphi}^{u}_\nu)$.

\

We now prove that determining a formula is prevalent. Here and in the
rest of the paper we often make use of sequences of pure extensions.
This gives us some closedness that the forcing as a whole does not have;
the situation is similar to that of Prikry forcing. We should think of
pure extensions as mild ones.

A \defn{pure sequence} is a sequence $\seq{\mu_i}_{i<\d}$ such that for
all $i<j < \d$, $\mu_j \pure \mu_i$. If $\d$ is limit, then such a
sequence has a natural limit (which by our notational conventions we
usually denote by $\mu_{<\d}$). For all $i<\d$ we have $\mu_{<\d} \pure
\mu_i$. However we note that it may be that $\mu_{<\d}$ is not a
condition in $\Q$ as its support may be too large. If $\d= \g+1$ then we
let $\mu_{<\d} = \mu_\g$.

\begin{lemma} Let $\mu_0\in \Q_u$ and let $\U\subset \Q_u$ be dense and
open. Then there is some $\mu\pure \mu_0$ in $\Q_u$ which determines
$\U$.
\label{lem: pure density lemma}
\end{lemma}

\begin{proof}
We construct a pure sequence $\seq{\mu_i}$, starting with $\mu_0$. If
$\mu_j$ is defined then we also pick some $a_j \in \B_{\mu_j}$ such that
$\mu_j\| a_j \in\U$ and for all $i<j$ we have $\mu_j(a_j\cap a_{i})=0$.
(Note that for all $i<j$, $\mu_j(a_i) = \mu_i(a_i)$.)

We keep constructing until we get stuck: we get some $\d$ such that
$\mu_{<\d}$ is defined but $\mu_{<\d}$ does not have any pure extension
$\vas$ such that there is some $a\in \B_\vas\cap \U$ which is
$\vas$-disjoint from all $a_i$ for $i<\d$.

We get stuck at a countable stage. For if not, $\seq{a_i}_{i<\w_1}$ are
pairwise $\mu_{<\w_1}$-disjoint which is impossible. This shows that at
limit stages $i$ we indeed have $\mu_{<i} \in \Q_u$ so the construction
can continue.

Suppose that we got stuck at stage $\d$; let $\mu = \mu_{<\d} \in \Q_u$.
We show that $\mu$ is as desired. $\{ \mu\| a_i\,:\, i<\d\}\subset \U$;
we claim that this is a maximal antichain in $\B_\mu$. If not, find some
$a\in \B_\mu$ which is $\mu$-disjoint from all $a_i$. Now there is some
extension of $\mu\|a$ in $\U$; it is of the form $\vas \|b$ where $\vas
\pure \mu$ and $b\subset a$. But then we can pick $\vas$ for $\mu_\d$
and $b$ for $a_\d$.
\end{proof}

\begin{lemma} Suppose that $\k<\d$ are inaccessible. Let $\mu_0\in \Q_{[\k,\d)}$ and
let $\U\subset \Q_\d$ be dense and open. Then there is some $\mu\pure
\mu_0$ in $\Q_{[\k,\d)}$ such that
\[ \{ \nu \in \Q_\k\,:\, \nu\mu \text{ determines }\U\} \]
is dense in $\Q_\k$.
\label{lem: dense determinacy}
\end{lemma}

\begin{proof}
We construct a pure sequence $\seq{\mu_i}$ of elements of $\Q_{[\k,\d)}$
of length below $\k^+$, starting with $\mu_0$. Together with this
sequence we enumerate an antichain $\A\subset \Q_\k$. At stage $i$, we
search for a pure extension $\vro$ of $\mu_{<i}$ in $\Q_\d$ which
determines $\U$ and is of the form $\vro = \nu'\mu'$ where $\nu' \in
\Q_\k$, $\mu' \in \Q_{[\k,\d)}$ and $\nu'$ is incompatible with all
elements enumerated so far into $\A$. If such exist, then we pick one,
enumerate $\nu'$ into $\A$ and let $\mu_{i} = \mu'$. If none such exist
then we stop the construction and let $\vas = \mu_{<i}$.

We must stop at some stage $i^*<\k^+$ because $|\Q_\k| = \k$.

Let $\nu \in \A$. If $\nu$ is enumerated into $\A$ at stage $i<i^*$ then
$\nu\mu_i$ determines $\U$; as $\vas \pure \mu_i$ we have $\nu\vas \pure
\nu\mu_i$ so $\nu\vas$ determines $\U$. It thus remains to show that
$\A$ is a maximal antichain of $\Q_\k$. Suppose not; let $\ups\in \Q_\k$
be incompatible with all elements of $\A$. By lemma \ref{lem: pure
density lemma}, we can find some $\vro \pure \ups\vas$ which determines
$\U$. We can write $\vro$ as $\nu'\mu'$ where $\nu'\pure \ups$ is in
$\Q_\k$ and $\mu'\pure \vas$ is in $\Q_{[\k,\d)}$. But $\nu'$ is
incompatible with all elements of $\A$ so we can pick $\mu_{i^*} =
\mu'$, which we didn't.
\end{proof}

\begin{scenario} \label{scn: the scenario}
Suppose now that $\k<\d$ are both inaccessible. Let $\bar \mu =
\seq{\mu_\a}_{\a<\a^*}$ be a pure sequence of measures in
$\Q_{[\k,\d)}$.

Let $G\subset \Q_\k$ be generic over $V$. For all $\nu\in G$, by lemma
\ref{lem: random embeds}, $G_\nu = G \cap \B_\nu$ is generic for
$\B_\nu$ over $V$. The system $\seq{G_\nu}_{\nu\in G}$ coheres: if
$\nu\le \vro$ then $G_\vro = {(i_{\vro}^{\nu})}^{-1} G_\nu$.

Let $\D_G = \{ \nu\mu_\a\,:\, \nu\in G \andd \a<\a^*\}$. This is a
directed system (under $\ge_\Q$). Note that from $\vas\in \D_G$ we can
recover $\nu$ and $\mu_\a$. We thus let, for $\vas=\nu\mu_\a\in \D_G$,
$\tau_\vas = \s_{\nu,\mu_\a}\circ \p_{\nu,\mu_\a} \colon \B_\vas \to
\B_{\mu_\a^{V[G_\nu]}}$ be the quotient by $G_\nu$ (this of course
depends on $G_\nu$ and not on $\vas$ alone, but we suppress its
mention). Lemmas \ref{lem: diagram 1}, \ref{lem: diagram 2} and
\ref{lem: diagram 3} and the discussion between them show that for any
$\vas=\nu\mu_\a\ge \vas'=\nu'\mu_\b$ in $\D_G$ and any $a\in \B_\vas$ we
have $\mu_\a^{V[G_\nu]}(\tau_\vas(a)) = \mu_\b^{V[G_{\nu'}]}
(\tau_{\vas'}(i_\vas^{\vas'}(a)))$.

\medskip

Let $\vphi$ be a formula of the forcing language for $\Q_\d$. For
$\vas=\nu\mu_\a\in \D_G$ we let $\xi_\vas(\vphi) = \mu_\a^{V[G_\nu]}
(\tau_\vas(\bool{\vphi}^\d_{\vas}))$. The analysis above shows that if
$\vas\ge \vas'$ are in $\D_G$ then $\xi_\vas(\vphi) \le
\xi_{\vas'}(\vphi)$, and that if $\vas$ $\d$-determines $\vphi$ then
$\xi_\vas(\vphi) = \xi_{\vas'}(\vphi)$ for all $\vas'\le \vas$. We
therefore let $\xi_G(\vphi) = \sup_{\vas\in \D_G} \xi_\vas(\vphi)$. To
calculate $\xi_G(\vphi)$ it is sufficient to take the supremum of
$\xi_\vas(\vphi)$ over a final segment of $\vas\in \D_G$ (or in fact any
cofinal subset of $\D_G$). If some $\vas\in \D_G$ determines $\vphi$
then $\xi_\vas(\vphi)$ is eventually constant and we get $\xi_G(\vphi) =
\max_{\vas\in \D_G} \xi_\vas(\vphi)$ which equals $\xi_\vas(\vphi)$ for
any $\vas$ which determines $\vphi$.
\end{scenario}

\begin{remark} This is important. Suppose that $M$ is an inner model of $V$. Then we
can work with this scenario ``mostly in $M$": we'll have all the
ingredients in $M$ (so $\k,\d$ are inaccessible in $M$, $\Q_\k, \Q_\d$
are in the sense of $M$) but the sequence $\bar \mu$ will not be in $M$.
Thus if $G\subset \Q_\k^M$ is generic \emph{over $V$} then the entire
system $(\D_G, \tau_\vas, \xi_\vas(\vphi),\dots)$ will be in $V[G]$ but
not in $M[G]$ (of course $G$ is generic over $M$ too). We can still
make, in $V[G]$, the above calculations of $\xi_G(\vphi)$ for $\vphi\in
M$ (although ``determining" and the calculation of $
\bool{\vphi}^\d_{\vas}$ and $\xi_\vas(\vphi)$ for each particular $\vas$
will be done in $M$ or $M[G]$).
\label{rmk: scenario in blue and black}
\end{remark}

\subsection{Real-valued measurability}

In this section we prove the following:

\begin{theorem} \label{thm: new forcing works}
Suppose that there is an elementary $j\colon V\to M$ with critical point
$\k$ such that $M^{<2^\k}\subset M$ (for example, if $\k$ is measurable
and $2^\k = \k^+$.) Then in $V^{\Q_\k}$, $\k$ is real-valued measurable.
\end{theorem}

Let $j$ be as in the theorem.

Let $\P = \Q_\k$. Then $\P\in M$ and $\P = \left(\Q_\k\right)^M$; more
importantly, $j(\P) = \left(\Q_{j(\k)}\right)^M$ (note that this is not
absolute; we do not have $j(\P) = \Q_{j(\k)}$). Let $\P' =
\left(\Q_{[\k,j(\k))} \right)^M$.

\medskip

What we do now is construct a pure sequence $\bar\mu=
\seq{\mu_\a}_{\a<2^\k}$ of elements of $\P'$. We start with a list
$\seq{\U_\a}_{\a<2^\k}$ of dense subsets of $j(\P)$ each of which is in
$M$ (note that this sequence is not in $M$). Rather than specify now
which dense sets we put on this list, we will, during the verifications
that $\k$ is real-valued measurable in $V^{\Q_\k}$, list dense sets that
are necessary for the proofs, making sure that we never put more than
$2^\k$ sets on the list.

Given $\seq{\U_\a}$, we construct $\bar \mu$ as follows. For $\mu_0$ we
pick any element of $\P'$. At stage $\a<2^\k$, we note that by the
closure property of $M$, $\seq{\mu_\b}_{\b<\a} \in M$ and so $\mu_{<\a}
\in M$. As $2^\k \le (2^\k)^M$ is less than the least inaccessible
beyond $\k$ in $M$, $\mu_{<\a}\in \P'$. We now apply lemma \ref{lem:
dense determinacy} in $M$, with $\k$ standing for $\k$, $j(\k)$ standing
for $\d$, $\mu_{<\a}$ for $\mu_0$ and $\U_\a$ for $\U$. The resulting
measure is $\mu_\a$.

If $G\subset \P$ is generic over $V$ then we find ourselves in scenario
\ref{scn: the scenario} (as modulated by remark \ref{rmk: scenario in
blue and black}). For every $\U$ on our list, we know that some $\vas\in
\D_G$ determines $\U$.

Let $\NN$ be the set of all $\P$-names for subsets of $\k$ (up to
equivalence); note that because $|\Q_\k|=\k$, $|\NN| = 2^\k$. If $G$ is
generic over $V$ then for every $A\in \NN$ we let $f_G(A) = \xi_G(\k\in
j(A))$. For the rest of this section, let $\d = j(\k)$.

\begin{lemma} Let $\nu\in \P$ and $A,B \in \NN$. Suppose that $\nu \force_\P
A\subset B$. Then $\nu \force_\P f_G(A) \le f_G(B)$. \end{lemma}
\begin{proof}
As we had in our discussion of Solovay's construction, $j\rest \P$ is
the identity. So in $M$, $\nu \force_{j(\P)} j(A) \subset j(B)$. Let
$G\subset \P$ be generic and suppose that $\nu\in G$. For any
$\vas=\nu'\mu_\a\in \D_G$ such that $\nu'\le \nu$ we have $\vas \le \nu$
so in $M$, $\vas \force_{j(\P)} j(A) \subset j(B)$ so $\bool{\k\in
j(A)}^\d_\vas \le_{\B_\vas} \bool{\k\in j(B)}^\d_\vas$ so
$\xi_\vas(\k\in j(A)) \le \xi_\vas(\k\in j(B))$. As this is true for a
final segment of $\vas\in \D_G$ we have $\xi_G(\k\in j(A)) \le
\xi_G(\k\in j(B))$. [Note that in this proof we didn't need any
particular $\U$.]
\end{proof}

It follows that $f_G$ induces a function on subsets of $\k$ in $V[G]$
(rather than only on their names). We show this function is the desired
measure on $\k$.

\begin{lemma} Let $A\subset \k$ be in $V$.
If $\k\in j(A)$ then $\force_\P f_G(A)=1$ and if $\k\notin j(A)$ then
$\force_\P f_G(A)=0$.
\end{lemma}

\begin{proof}
Suppose that $\k\in j(A)$. Then in $M$, every condition in $j(\P)$
forces this fact. Let $G\subset \P$ be generic. It follows that for all
$\vas\in \D_G$, $\bool{\k\in j(A)}^\d_\vas = 1_{\B_\vas}$ so
$\xi_G(\k\in j(A)) =1$.

We get a similar argument if $\k\notin j(A)$.
\end{proof}

\begin{lemma} Let $\seq{B_n}_{n<\w}$ be a sequence of names in $\NN$.
Suppose that $\nu \in \P$ forces that $B_{n}$ are pairwise disjoint.
Then $\nu\forces_\P f_G(\cup_n B_n) = \sum_{n<\w} f_G(B_n)$.
\end{lemma}

\begin{proof}
Let $B\in \NN$ be such that $\nu\forces_\P B = \cup_n B_n$.

We have $j(\nu)=\nu$ and $j(\seq{B_n}_{n<\w}) = \seq{j(B_n)}_{n<\w}$; so
in $M$, $\nu$ forces (in $j(\P)$) that $j(B) = \cup_{n<\w} j(B_n)$ is a
disjoint union. Again let $G$ be generic such that $\nu\in G$.

Let $\U$ be the collection of $\mu\in j(\P)$ extending $\nu$ such that
in $M$, if $\mu\force_{j(\P)} \k\in j(B)$ then for some $n<\w$,
$\mu\force_{j(\P)} \k\in j(B_n)$. $\U\in M$ and $\U$ is dense in
$j(\P)$. We assume that some $\vas\in \D_G$ (and so a final segment of
$\vas\in \D_G$) determines $\U$. [Note that the number of such sequences
$\seq{B_n}$ is $|\NN|^{\aleph_0} = 2^\k$ so we may put all the
associated $\U$'s on the list.]

If $\vas$ determines $\U$ then $\bool{\k\in j(B)}^\d_\vas =
\sum^{\B_\vas}_{n<\w} \bool{\k\in j(B_n)}^\d_\vas$. Also, if $\vas\in
\D_G$ extends $\nu$ then for $n\ne m$ we have $\bool{\k\in
j(B_n)}^\d_\vas \meet_{\B_\vas} \bool{\k\in j(B_m)}^\d_\vas = 0$. It
follows that in addition, if $\vas$ determines $\k\in j(B)$ then it
determines $\k\in j(B_n)$ for every $n<\w$ (again only $2^\k$ many
$\U$'s to add).

Thus, for plenty $\vas\in \D_G$ we have $f_G(B) = \xi_\vas(\k\in j(B)) =
\sum_{n<\w} \xi_\vas(\k\in j(B_n)) = \sum_{n<\w} f_G(B_n)$ as required.
\end{proof}

\begin{lemma} Suppose that $\g<\k$ and that $\seq{B_\a}_{\a<\g}$ is a sequence of
names in $\NN$.  Suppose that $\nu\in \P$ forces that for all $\a<\g$,
$f_G(B_\a)=0$. Then $\nu\forces_\P f_G(\cup_\a B_\a)=0$.
\end{lemma}
\begin{proof}
Let $B\in \NN$ be such that $\nu\forces_\P B = \cup_{\a<\g} B_\a$. Then
in $M$, $\nu\forces_{j(\P)} j(B) = \cup_{\a<\g} j(B_\a)$. Let $G\subset
\P$ be generic over $V$ and suppose that $\nu\in G$. For all $\a<\g$,
$f_G(B_\a) = 0$ so for all $\vas\in \D_G$, $\bool{\k\in j(B_\a)}^\d_\vas
= 0_{\B_\vas}$.

Let $\U$ be the collection of conditions $\mu\in j(\P)$ extending $\nu$
such that in $M$, if $\mu\forces_{j(\P)} \k\in j(B)$ then for some
$\a<\g$, $\mu\forces_{j(\P)} \k \in j(B_\a)$. Then $\U\in M$ and $\U$ is
dense below $\nu$ in $j(\P)$. We assume that some $\vas\in \D_G$
determines $\U$. If $\vas$ determines $\U$ then $\bool{\k\in
j(B)}^\d_\vas = \sum^{\B_\vas}_{\a<\g} \bool{\k\in j(B_\a)}^\d_\vas =
0$. Thus on a final segment of $\vas\in \D_G$ we have $\bool{\k\in
j(B)}^\d_\vas=0$ so $\xi_G(\k\in j(B))=0$.

Now there are $|\NN|^{<\k}= 2^\k$ such sequences $\seq{B_\a}$ so we only
need $2^\k$ many such $\U$ on our list of dense sets to determine.
\end{proof}

%%%%%%%%%%%%%%%%%%%%%%%%%%%%%%%%%%%%%%%%%%%%%%%%%%%%%%%%%%%%%%%%%%%%%%%%%%%

\section{General sequences}

To facilitate the definition, we introduce some notation. Suppose that
$w\subset \On$ and that $\bar x = \seq{x_\a}_{\a\in w}$ is a sequence of
reals. Suppose that $B\subset \R^{\otp w}$. Then we say that $\bar x \in
B$ if $\seq{x_{f(\xi)}}\in B$, where $f\colon \otp w \to w$ is
order-preserving. If $\bar B = \seq{B_i}_{i<\s}$ is a sequence of sets
such that for all $i<\s$, $B_i \subset \R^i$, then we say that $\bar
x\in \bar B$ if $\bar x \in B_{\otp w}$.

Let $\s\le \k$ be regular, uncountable cardinals.

\begin{definition} A \defn{$\k$-null}
set is a union of fewer than $\k$ null sets. \footnote{Let $\mathcal{N}$
be the ideal of null sets. If $\kappa$ is real-valued
       measurable, we have $non(\mathcal{N})=\aleph_1$ and
       $cov(\mathcal{N})\ge\kappa$ (\cite{Fremlin}). Hence, for a real-valued
       measurable $\kappa$, $\kappa$-null sets form a proper ideal
       extending $\mathcal{N}$ properly. By the inequality above we have
       $cov(\mathcal{N})=\kappa$ in Solovay's model as well as in the new model.
       The existence of a $\s$-general sequence, which separates between the models,
       can be viewed as a strengthening of the equation $cov(\mathcal{N})=\kappa$.}
\end{definition}

\begin{definition} A \defn{$\k$-null sequence} is a sequence
$\bar B = \seq{B_i}_{i<\s}$ such that for each $i<\s$, $B_i$ is a
$\k$-null subset of $\R^i$.
\end{definition}

\begin{definition} Let $A$ be any set.
A \defn{noncountable club} on $[A]^{<\s}$ is some $\CC\subset [A]^{<\s}$
which is cofinal in $\left([A]^{<\s},\subseteq\right)$ and is closed
under taking unions of increasing chains of uncountable cofinality.
\end{definition}

\begin{definition} Let $\bar B = \seq{B_i}_{i<\s}$ be such that for
all $i<\s$, $B_i\subset \R^i$. Let $\bar x = \seq{x_\a}_{\a\in U}$ be a
sequence of reals (where $U\subset \On$). We say that $\bar x$
\defn{escapes} $\bar B$ if there is some noncountable club $\CC$ on
$[U]^{<\s}$ such that for all $w\in \CC$, $\bar x\rest w \notin \bar B$.
\end{definition}

\begin{definition} \label{def: general sequence}
A sequence $\bar x = \seq{x_\a}_{\a<\k}$ of reals is \defn{$\s$-general}
if for every $\k$-null sequence $\bar B = \seq{B_i}_{i<\s}$, there is
some final segment $W$ of $\k$ such that $\bar x \rest W$ escapes $\bar
B$.
\end{definition}

\subsubsection{Justifying the definition}

Na\"{i}ve approaches might have liked to \linebreak strengthen the above
definition. However, it is fairly straightforward to see that expected
modes of strengthening result in empty notions. For example, one would
like to eliminate the restriction to a final segment of $\k$. But given
a sequence $\bar x= \seq{x_\a}_{\a<\k}$, we can let, for $i<\s$, \[B_i =
\{ x_{0}\} \times \R^{i\setminus\{0\}} = \{ \bar y \in \R^i\,:\, \bar
y(0) = x_{0}\}.\] Then whenever $w\subset \k$ such that $0\in w$, $\bar
x\rest w \in B_{\otp w}$ (and every noncountable club on $[\k]^{<\s}$
contains such a $w$.)

Accepting the restriction to a final segment, we may ask why we need to
restrict to a club - why we can't have $\bar x\rest w\notin \bar B$ for
all $w\in [W]^{<\s}$. But consider \[B_\w = \{\bar y\in \R^\w\,:\,
\forall n<\w \,\,\, \bar y(2n)(0) = \bar y (2n+1)(0)\}.\] Given a final
segment $W$ of $\k$, we can always choose some $\w$-sequence $w\subset
W$ such that $\bar x\rest w\in B_\w$.

\subsection{General sequences in Solovay's model}

\begin{theorem} Let $\k$ be inaccessible. Then
in $V^{\B_{m_\k}}$, the random sequence is $\s$-general for all regular,
uncountable $\s<\k$.
\label{thm: the random sequence is general}
\end{theorem}

This relies on the following well-known fact:
\begin{fact} Let $\P$ be a notion of forcing which has the $\l$-Knaster
condition for all regular uncountable $\l<\s$, and let $A\in V$. Then
(in $V^\P$), $\left([A]^{<\s}\right)^V$ is a noncountable club of
$[A]^{<\s}$.
\end{fact}

(Recall that $\P$ has the $\l$-Knaster condition if for all $A\subset
\P$ of size $\l$, there is some $B\subset A$ of size $\l$ such that all
elements of $B$ are pairwise compatible in $\P$.)

\begin{proof}
Let $\A = \left([A]^{<\s}\right)^V$. To see that $\A$ is cofinal in
$[A]^{<\s}$, let $u$ be a name for an element of $[A]^{<\s}$. Let $p\in
\P$ force that $\{t_i\,:\, i<\l\}$ is an enumeration of $u$ (for some
$\l<\s$). For $i<\l$, let $P_i\subset \P(\le p)$ be a maximal antichain
of elements $q$ which force that $t_i = a_{i,q}$ for some $a_{i,q}\in
A$. Then $p$ forces that $w = \{a_{i,q}\,:\, i<\l, q\in P_i\}$ (which is
in $\A$) contains $u$.

Now suppose that $p\in \P$ forces that $\seq{u_i}_{i<\l}$ is an
increasing sequence in $\A$, for some regular uncountable $\l<\s$. For
every $i<\l$ pick some $p_i\le p$ and some $w_i\in \A$ such that
$p_i\forces u_i=w_i$. Note that if $p_i$ and $p_j$ are compatible and
$i<j$ then $w_i \subset w_j$. Let $X\in [\l]^\l$ be such that for
$i,j\in X$, $p_i$ and $p_j$ are compatible. Without loss of generality,
assume that $\P$ is a complete Boolean algebra. For $i\in X$ let $q_i =
\sum^\P_{j>i, j\in X} p_j$. Then $\seq{q_i}_{i\in X}$ is decreasing and
so halts at some $q_{i^*}$. Then $q_{i^*}$ forces that for unboundedly
many $i\in X$, $p_i\in G$, and so that $\cup_{i<\l} u_i = \cup_{i<\l}
w_i$ which is in $\A$.
\end{proof}

\begin{fact} A measure algebra $(\B,\mu)$ has the $\l$-Knaster condition
for all regular uncountable $\l$. \end{fact}

\begin{proof}
This is well-known; see, for example, \cite{Argyros_Kalamidas}. We give
a proof for the sake of completeness.

Suppose that $\{b_i\,:\, i<\l\}\subset \B$. Let $X_0\in [\l]^\l$ such
that for all $i\in X_0$, $\mu(b_i)>1/n$. Inductively define $X_{m+1}$
from $X_m$: if there is some $i\in X_m$ such that for $\l$ many $j\in
X_m$, $b_i\cap b_j =0$, then let $i$ be minimal such and let $X_{m+1} =
\{j\in X_m\,:\, j>i \,\,\,\&\,\,\,b_j\cap b_i =0\}$. This process has to
terminate with some $X_{m^*}$ because $\sum_{i\in X_m} b_i - \sum_{i\in
X_{m+1}} b_i$ has measure $>1/n$. We can now find $Y\in [X_{m^*}]^\l$
which indexes a set of pairwise compatible conditions by inductively
winnowing all $j$ such that $b_j$ is disjoint from something we put into
$Y$ so far.
\end{proof}

\begin{proof}[Proof of Theorem \ref{thm: the random sequence is general}]
Let $G\subset \B_{m_\k}$ be generic over $V$, and let $\bar r =
\seq{r_\a}_{\a<\k}$ be the random sequence obtained from $G$. In $V[G]$,
let $\bar B = \seq{B_i}_{i<\s}$ be a $\k$-null sequence of length $\s$.
For $i<\s$ choose null sets $B^i_{\a}\subset \R^i$ for $\a<\a_i<\k$ such
that $B_i = \cup_{\a<\a_i} B^i_\a$.

A code for each $B^i_\a$ is a real, together with some countable subset
of $i$. It follows that there is some $\theta<\k$ such that each
$B^i_\a$ is defined in $V'=V[G\cap \B_{m_\theta}]$. Let $W =
[\theta,\k)$. Then $\bar r\rest W$ is random (for $\B_{m_W}$) over $V'$.
Let $w\in \left([W]^{<\s}\right)^{V'}$, and let $i = \otp w$. The
collapse $h\colon w\to i$ induces a bijection $h\colon \R^w \to \R^i$.
Let $\a<\a_i$ and consider $h^{-1} B^i_{\a}$; this is a null subset of
$\R^w$ defined in $V'$ and so $\bar r\rest w$, being random over $V'$,
is not in $h^{-1} B^i_\a$. Which means, in our notation, that $\bar r
\rest w\notin B^i_\a$ and so $\bar r\rest w \notin \bar B$. The
noncountable club $\CC = \left([W]^{<\s}\right)^{V'}$ thus witnesses
that $\bar r\rest W$ escapes $\bar B$.
\end{proof}

\subsection{Some necessary facts about $\Q_\k$}

The following information will be useful in showing the lack of general
sequences. From now, assume that $\aleph_2\le \s< \k$, both $\s$ and
$\k$ are regular, and that $\s$ is at most the least inaccessible; this
is a convenience, since then $\Q_\k$ is purely $\s$-closed.

\subsubsection{Cardinal preservation}

\begin{lemma} All cardinals and cofinalities below the least inaccessible are preserved
by $\Q_\k$. \end{lemma}

This is important; if $\s$ is not regular in the extension then
$[W]^{<\s}$ ceases to be interesting.

\begin{proof}
Let $\theta$ be a regular, uncountable cardinal below the least
inaccessible cardinal. Let $\l<\theta$ and suppose that $\mu_0\in \Q_\k$
forces that $f\colon \l\to \theta$ is a function. By lemma \ref{lem:
pure density lemma} construct a pure sequence $\seq{\mu_i}_{i<\l}$ in
$\Q_\k$ starting with $\mu_0$ such that for each $i<\l$, $\mu_i$
determines the value of $f(i)$ (that is, the collection of $a\in
\B_{\mu_i}$ such that for some $\g<\theta$, $\mu_i \| a \force_{\Q_\k}
f(i) = \g$ is dense in $\B_{\mu_i}$). For $i<\l$ let $A_i$ be the
(countable) set of such values $\g$. For every $i$, $\B_{\mu_i}
\compembed \Q_\k (\le \mu_{i})$ and so $\mu_{<\l}$ forces that the range
of $f$ is contained in $\cup_{i<\l} A_i$.
\end{proof}

As for preservation of cardinals beyond the least inaccessible, we
mention that the only cardinals that may be collapsed by $\Q_\k$ are
those that lie between $\d$ and $2^\d$, where $\d<\k$ is a singular
limit of inaccessible cardinals. We omit the proof as it does not
involve new techniques (and we do not use this fact). As to whether
$\Q_\k$ does collapse any cardinals, it seems this may be independent.
The general results of \cite{DjSh:659} may be relevant here.

\subsubsection{Finding elements of clubs}

\begin{lemma} Let $A\in V$. Suppose that $\mu\in \Q_\k$ forces that
$\CC$ is a noncountable club on $[A]^{<\s}$. Then there is some (pure)
extension $\nu$ of $\mu$ and some $w\in [A]^{<\s}$ (in $V$) such that
$\nu\forces w\in \CC$.
\label{lem: getting at clubs}
\end{lemma}

\begin{proof}
We show the following claim: given $\mu\in \Q_\k$ forcing that $\CC$ is
a noncountable club on $[A]^{<\s}$ and given some $w\in [A]^{<\s}$,
there is some $\nu$ purely extending $\mu$ and some $w'\in [A]^{<\s}$
containing $w$ such that $\nu$ forces that there is some $v\in \CC$,
$w\subset v \subset w'$.

This suffices: given $\mu$ as in the lemma, we construct a pure sequence
$\seq{\mu_\a}_{\a<\w_1}$ starting with $\mu$ and an increasing sequence
of $w_\a \in [A]^{<\s}$ such that $w_0 =0$, and $\mu_{\a}$ forces that
there is some $v_\a\in \CC$, $w_{<\a} \subset v_\a \subset w_{\a}$. Then
$\mu_{<\w_1}$ forces that $\cup_{\a<\w_1} w_\a = \cup_{\a<\w_1} v_\a$ is
in $\CC$.

So let $\mu$, $w$ be as in the claim. Let $w^*$ be a name such that
$\mu\forces w^*\in \CC$ and $w\subset w^*$ ($\CC$ is cofinal). First,
let $\mu'$ be a pure extension of $\mu$ such that there is an antichain
$\seq{a_n}_{n<\w}$ of $\B_{\mu'}$ and cardinals $\l_n <\s$ such that
$\mu'\|a_n \forces |w^*| = \l_n$; for every $n$, let
$\seq{x^n_i}_{i<\l_n}$ be a list of names such that $\mu'\|a_n \forces
w^* = \{x^n_i\,:\, i<\l_n\}$. We now construct a pure sequence $\mu_i$
for $i< \l = \sup_n \l_n$; for each $i<\l$ and each $n$ such that
$i<\l_n$, the collection of $b\in \B_{\mu_i}$ such that $b\subset a_n$
and for some $a\in A$, $\mu_i\|b \forces a = x^n_i$ is dense below
$a_n$; there are only countably many such $a$. Then $\mu_{<\l}$ forces
that $w^*$ is contained in $w'$, the collection of all such $a$'s which
appeared in the construction ($\s$ is regular, so $\l<\s$ and so $w'$
has size $<\s$).
\end{proof}

In fact, for every $w\in [A]^{<\s}$ and such $\mu$, there is a pure
extension forcing that some $w'$ containing $w$ is in $\CC$. This is
immediate from the proof, or from the fact that $\{v\in \CC\,:\,
v\supset w\}$ is also a noncountable club of $[A]^{<\s}$ (in the
extension).

\subsubsection{Approximating measures by pure measures}

The following is an easy fact which follows from regularity of our
measures:

\begin{lemma} Let $\mu$ be a pure measure, and let $B\in \B_{\mu}$. Then
for all $\e<1$ there is some \emph{pure} measure $\nu$ which is a
localization of $\mu$  such that $\nu(B)>\e$.
\label{lem: thickness lemma}
\end{lemma}

\begin{proof}
This follows from regularity of $\mu$. There is some open set $U\supset
B$ such that $\mu(U\setminus B) < (1-\e)\mu(B)/\e$; so
$\mu(B)/\mu(U)>\e$. We can present $U$ as a disjoint union of cylinders
$U_n$; for some $n$ we must have $\mu(U_n\cap B)/\mu(U_n) >\e$. Then
$\mu\| U_n$ is a pure measure and is as required.
\end{proof}

We need a certain degree of uniformity.

\begin{lemma} Let $B \subset \R^2$ be a positive Borel set. Then there is
some positive $A\subset \R$ such that for all $\e<1$ there is some
positive $C\subset \R$ such that $m_2\|(A\times C)(B) >\e$.
\label{lem: funny denisty}
\end{lemma}

\begin{proof}
Let $X\elementary V$ \footnote{Yes, we mean $X\elementary H(\chi)$.
Complaints are to be lodged with set models of ZFC.} be countable such
that $B\in X$. Let $C_0$ be the measure-theoretic projection of $B$ onto
the $y$-axis (of course $C_0\in X$). $C_0$ is positive, so we can pick
some $r^*\in C_0$ which is random over $X$. Let $A = B^{r^*} = \{x\in
\R\,:\, (x,r^*)\in B\}$ be the section defined by $r^*$; since $r^*\in
C_0$, $A$ is positive. Note that in $X$ there is a name for $B^{r^*}$,
where $r^*$ is a name for the generic random real.

Let $\d>0$ be in $X$. By regularity of Lebesgue measure, there is some
clopen set $U\subset \R$ such that $m(U \triangle A) < \d$.\footnote{Let
$V\supset A$ be open such that $m(V\setminus A)< \d/2$; and recall that
every open set is an increasing union of clopen sets.} Of course $U\in
X$. Then there is some positive $C\subset C_0$ in $X$ such that
$C\forces_{\B_m} m(U\triangle B^{r^*}) < \d$.

For almost all $r\in C$ (those that are random over $X$), we have
$m(U\triangle B^{r})<\d$. For such $r$, $m(A\setminus B^r) \le
m(A\setminus U) + m(U\setminus B^r) \le 2\d$. So by Fubini's theorem,
$m_2(A\times C \setminus B) = \int_C m(A\setminus B^r)\,d\,r \le 2\d
m(C)$; we get that $m_2\|(A\times C)(\lnot B) \le 2\d/m(A)$.\footnote{We
glossed over uses of the forcing theorem over $X$, which is not
transitive. We really work with $X$'s collapse and use absoluteness. For
example, we got $C\in X$ such that $C\forces_{\B_m} m(U\triangle
B^{r^*})<\d$. Let $\p\colon X\to M$ be $X$'s transitive collapse. Then
in $M$, $\p(C)$ forces (in $\B_m^M$) that $m(\p(U)\triangle
\p(B)^{r^*})<\d$. If $r\in C$ is random over $M$, then in $M[r]$,
$m(\p(U)\triangle \p(B)^r)<\d$. But $\p(B)^{M[G]} = B \cap M[G]$ and
similarly for $U$. Thus indeed $m(U\triangle B^r)<\d$ as we claimed.}
\end{proof}

\begin{corollary} Let $\vro,\mu\in \Q_\k$ and let $\mu$ be pure; assume $u^\vro\cap u^\mu = 0$.
Let $B\in \B_{\vro\mu}$. Then there is a localization $\vro'$ of $\vro$
such that for all $\e<1$ there is some pure $\mu'$ which is a
localization of $\mu$ and such that $\vro'\mu'(B) > \e$.
\label{cor: funny density}
\end{corollary}

\begin{proof} What we need to note is that the proof of the previous
lemma holds for $\vro \times \mu$ (in place of $m \times m$) (we just
use the relevant measure algebra); we get a set $A\in \B_\vro$ such that
for all $\d>0$ there is some $C\in \B_\mu$ such that $\vro\mu\|(A\times
C)(B) > 1-\d$.

Fix some $\d>0$. Get the appropriate $C$; we have
\[ \frac{\vro\mu(A\times C \setminus B)}{\vro\mu(A\times C)} < \d. \]
By Lemma \ref{lem: thickness lemma}, we can find some cylinder $\tilde C
\in \B_\mu$ sufficiently close to $C$ so that both $\mu \| \tilde C (C)
> 1-\d$ and $\mu(C) / \mu(\tilde C) < 1+\d$; from the first we get
\[ \frac{\mu(\tilde C \setminus C)}{\mu(\tilde C)} < \d .\]
Note that $A\times \tilde C \setminus B \subset A\times (\tilde C
\setminus C) \cup (A\times C \setminus B)$. Combining everything, we get
\begin{align*}
 \frac{\vro\mu(A\times \tilde C \setminus B)}{\vro\mu(A\times \tilde C)} \,\le\,
    \frac{\vro\mu(A\times (\tilde C \setminus C))\, +\, \vro\mu(A\times
                    C\setminus B)}{\vro\mu (A\times \tilde C)} \, = \,
                    \\
    \frac{\vro(A)\mu(\tilde C\setminus C)}{\vro(A)\mu(\tilde C)} +
            \frac{\vro\mu(A\times C \setminus B)}{\vro\mu(A\times C)}
                \cdot \frac{\vro(A)\mu(C)}{\vro(A)\mu(\tilde C)} \, \le \,
    \d + \d(1+\d) .
\end{align*}
We can thus let $\vro' = \vro\| A$ and $\mu' = \mu\|\tilde C$; the
latter is pure because $\tilde C$ is a cylinder. We get $\vro'\mu'(B)
\ge 1 - 2\d - \d^2$ which we can make sufficiently close to 1.
\end{proof}

\subsection{In the new model}

\begin{theorem} Suppose that $\k$ is Mahlo
for inaccessible cardinals, and that $\s\ge \aleph_2$ is at most the
least inaccessible (and is regular). Then in $V^{\Q_\k}$, there are no
$\s$-general sequences.
\label{thm: no general sequences}
\end{theorem}

In fact, we prove something stronger:

\begin{theorem} Suppose that $\k$ is Mahlo for inaccessible cardinals,
and that $\s\ge \aleph_2$ is regular, and is at most the least
inaccessible. Then there is a $\k$-null sequence $\bar B$ of length $\s$
such that in $V^{\Q_\k}$, no $\k$-sequence of reals $\bar r$ escapes
$\bar B$.
\label{thm: stronger no general sequences}
\end{theorem}

(We have here identified $\bar B$ as it is interpreted in $V$ and in
$V^{\Q_\k}$. Of course, for every $\k$-null $B$, if $B = \cup_{i<i^*}
B_i$ for some $i^*<\k$ then for any $W\supset V$ we let $B^{W} =
\cup_{i<i^*} B_i^W$.)

\begin{proof}
Work in $V$. We define $\bar B$ as follows: for $i<\s$, an increasing
$\w$-sequence $\bar j=\seq{j_n}_{n<\w}$ from $i$, and $k<2$, we let
\[ B^i_{\bar j,k} = \cap_{n<\w} [\seq{k}]^{j_n} = \{ \bar x\in \R^i\,:\,
    \forall n<\w\,\,\, (\bar x(j_n)(0)=k)\}; \]
and we let $B_i$ be the union of the $B^i_{\bar j,k}$ for all increasing
$\bar j$ from $i$ and $k<2$. Each $B^i_{\bar j,k}$ is null, and $\k$ is
inaccessible, so $B_i$ is $\k$-null. As $\k$ remains a cardinal in
$V^{\Q_\k}$, $B_i^{V^{\Q_\k}}$ is also $\k$-null in $V^{\Q_\k}$.

\medskip

Let $\mu^*\in \Q_\k$ force that $\bar r = \seq{r_\a}_{\a<\k}$ is a
sequence of reals and that $\CC$ is a noncountable club on $[\k]^{<\s}$.

For every $\g<\k$, find some $\mu_\g$ extending $\mu^*$ and some
$k(\g)\in 2$ such that $\mu_\g\forces r_\g(0) = k(\g)$.

Suppose that $\g> \sup u^{\mu^*}$. Then we can find $\vpi_\g\in \Q_\g$
which is an extension of $\mu^*$, a pure measure $\nu_\g\in
\Q_{[\g,\k)}$, and some Borel $B_\g$, such that $\mu_\g =
\left(\vpi_\g\nu_\g\right) \| B_\g$. Let $u_\g\subset u^{\mu_\g}$ be a
countable support for $B_\g$.

We now winnow the collection of $\mu_\g$'s. Let $S_0$ be the set of
inaccessible cardinals below $\k$ (but greater than $\sup u^{\mu^*}$);
for $\g\in S_0$ we have $\Q_\g\subset V_\g$. Thus on some stationary
$S_1\subset S_0$, the function $\g\mapsto \vpi_\g$ is constant. Next, we
find $S_2\subset S_1$ such that on $S_2$:
\begin{itemize}
\item $u_\g\cap \g$ and $\otp u_\g$ are constant;
\item Under the identification of one $\R^{u_\g\setminus \g}$ to the other by the
order-preserving map, $\nu_\g\rest (u_\g\setminus \g)$ is constant;
\item Under the identification of one $\R^{u_\g}$ to the other by the
order-preserving map, $B_\g$ is constant;
\item $k(\g)$ is a constant $k^*$.
\end{itemize}

By these constants, and using Corollary \ref{cor: funny density}, we can
find some $\mu^{**}$, a localization of $\vpi_\g$ for $\g\in S_2$, such
that for all $\e<1$ and all $\g\in S_2$, there is some pure $\vas$ which
is a localization of $\nu_\g$, such that $(\mu^{**}\vas) (B_\g) > \e$.

\medskip

We now amalgamate countably many $\mu_\g$'s in the following way. Pick
an increasing sequence $\seq{\g_n}_{n<\w}$ from $S_2$. For each $n<\w$,
let $\vas_n$ be a pure measure, which is a localization of $\nu_{\g_n}$,
such that $\mu^{**}\vas_n (B_{\g_n}) > q_n$, where $\seq{q_n}$ is a
sequence of rational numbers in $(0,1)$ chosen so that $\sum_{n<\w}
(1-q_n) < 1$. We note that $u^{\vas_n}$ are pairwise disjoint, and so we
can take their product $\vas^* = \mu^{**}\vas_0\vas_1\dots$. We now let
$\vas^{**} = \vas^{*} \| \cap_{n<\w} B_{\g_n}$; the $\vas_n$ were chosen
so that this is indeed a measure.

The point is that for all $n$, $\vas^{**} \le (\mu^{**}\vas_n)\|B_{\g_n}
\le (\mu^{**}\nu_{\g_n})\| B_{\g_n} \le \mu_{\g_n}$. Thus for all $n$,
$\vas^{**}\forces r_{\g_n}(0)=k^*$.

Finally, by Lemma \ref{lem: getting at clubs}, let $\vro$ be some
extension of $\vas^{**}$ which forces that some $w\in \CC$, where $w\in
V$ and $w\supset \{\g_n\,:\, n<\w\}$. Let $i = \otp w$ and let $h\colon
w\to i$ be the collapse. Define $\bar j$ by letting $j_n = h(\g_n)$.
Then $\vro$ forces that $\bar r\rest w \in B^i_{\bar j, k^*}$ so that
$\bar r\rest w\in \bar B$. Thus $\mu^*$ could not have forced that $\CC$
witnesses that $\bar r$ escapes $\bar B$.
\end{proof}

%%%%%%%%%%%%%%%%%%%%%%%%%%%%%%%%%%%%%%%%%%%%%%%%%%%%%%%%%%%%%%%%%%%

\bibliography{RVMbib,listb}
\bibliographystyle{elsart-num}

\end{document}